\documentclass[12pt]{amsart}
\usepackage{amscd,amssymb,amsmath}
\newtheorem{thm}{Theorem}[section]
\newtheorem{corol}[thm]{Corollary} 
\newtheorem{lemma}[thm]{Lemma}
\newtheorem{prop}[thm]{Proposition}
\newtheorem{ques}[thm]{Question}
\newtheorem{mainth}[thm]{Main Theorem}
\newtheorem{problem}[thm]{Problem}

\theoremstyle{definition}
\newtheorem{defin}[thm]{Definition}

\theoremstyle{remark}

\numberwithin{equation}{section}

\def\obr{^{-1}}
\def\sbs{\subset}

\def\Homeo{\mathrm{Homeo\,}}
\def\diam{\mathrm{diam\,}}
\def\Inn{\mathrm{Inn\,}}
\def\e{\epsilon}
\def\o{\omega}
\def\G{\Gamma}
\def\La{\Lambda}
\def\bull{\bullet}

\def\a{\alpha}
\def\d{\delta}
\def\D{\Delta}
\def\f{\varphi}
\def\g{\gamma}
\def\k{\kappa}

\def\Th{\Theta}

\def\Z{{\mathbf Z}}

\def\cl#1{\overline {#1}}

\def\Cal{\mathcal}
\def\sB{{\Cal B}}

\def\sF{{\Cal F}}
\def\sG{{\Cal G}}

\def\sL{{\Cal L}}
\def\sM{{\Cal M}}
\def\sN{{\Cal N}}

\def\sR{{\Cal R}}

\def\sU{{\Cal U}}

\def\ti{\times}

\def\Iso{\mathrm{Iso\,}}
\def\Is{\mathrm{Iso\,}}

\def\U{{\Bbb U}}

\begin{document}
\author[V. Uspenskij]{Vladimir V. Uspenskij}

\date{Version of 13 March 2007}

\address 
{Department of Mathematics, 321 Morton Hall, Ohio University, 
Athens, Ohio 45701, USA}

\email 
{uspensk@math.ohiou.edu}

\thanks{%
{\it 2000 Mathematics Subject Classification:}
Primary 22A05. Secondary 06F05, 22A15, 54D35, 
54E15, 54E50, 54H11, 54H12, 54H15, 57S05}

\keywords{Topological group, uniformity, semigroup, idempotent,
isometry, Urysohn metric space, Roelcke compactification, unitary group}

\title{On subgroups of minimal topological groups}

\begin{abstract}
A topological group is {\it minimal\/} 
if it does not admit a strictly coarser
Hausdorff group topology.
The {\it Roelcke uniformity\/} 
(or {\it lower uniformity\/}) on a topological group 
is the greatest lower bound
of the left and right uniformities.  A group is
{\it Roelcke-precompact\/}
if it is precompact with respect to the Roelcke uniformity.
Many naturally arising non-Abelian topological groups are Roelcke-precompact
and hence have a natural compactification. We use such compactifications
to prove that some groups of isometries are minimal.
In particular, if $\U_1$ is the Urysohn universal
metric space of diameter $1$, 
the group $\Iso(\U_1)$ of all self-isometries of $\U_1$
is Roelcke-precompact, topologically simple
and minimal.
We also show that every topological group is a subgroup
of a minimal topologically simple Roelcke-precompact group of the form
$\Iso(M)$, where $M$ is an appropriate non-separable version of
the Urysohn space.

The paper is available online: arXiv:math.GN/0004119.
\end{abstract}


\maketitle

\section{Introduction} \label{s:intro}

This paper was motivated by the following questions:

\begin{ques}[V. Pestov, A. Arhangelskii, 1980's]
\label{q:PesArh}
What are subgroups of minimal topological groups?
\end{ques}

\begin{ques}[W. Roelcke, 1990]
\label{q:roel}
What are subgroups of lower precompact topological groups?
\end{ques}

We now explain and discuss the notions of a minimal group
and of a lower precompact group.

Compact spaces $X$ can be characterized among all Tikhonov spaces by each of
the following two properties:
(1) $X$ is minimal, in the sense that $X$ admits no
strictly coarser Tikhonov (or Hausdorff) topology;
(2) $X$ is absolutely closed,
which means that  $X$ is closed in any Tikhonov space $Y$ containing $X$ as a
subspace.
One can consider the notions of minimality and absolute closedness also
for other classes of spaces. For example, for the class of Hausdorff spaces one
gets the notions of $H$-minimal and $H$-closed spaces which are no longer
equivalent to each other or to compactness but are closely related: a space is
$H$-minimal iff it is $H$-closed and semiregular,
and a space is compact iff it is
$H$-minimal and satisfies the Urysohn separation axiom. 
See the survey \cite{PorSte} for a discussion of these notions.

Let us now consider the case of topological groups. All topological groups  are
assumed  to be Hausdorff, unless otherwise explicitly stated.
A topological group is {\it minimal\/} if it does not
admit a strictly coarser Hausdorff group topology%
\footnote{The survey \cite{DikSurv} on minimal groups
contains a lot of information and more than a hundred references.}.
A topological group is  {\it
absolutely  closed\/}  if it is closed in every topological group containing it
as a topological subgroup. A topological group $G$ is absolutely closed if  and
only  if  it  is  {\it Rajkov-complete\/}, 
or {\it upper complete}, that is complete with respect to the
upper uniformity which is defined as the least upper bound $\sL\vee\sR$ of  the
left and the right uniformities 
on $G$. Recall that the sets $\{(x,y): x\obr y\in U\}$,
where  $U$  runs  over  a  base  at  unity  of $G$, constitute a base of
entourages for the left uniformity $\sL$ on $G$.  In  the  case  of  the  right
uniformity  $\sR$,  the condition $x\obr y\in U$ is replaced by $yx\obr \in U$.
We shall call Rajkov-complete groups simply {\it complete}.
The {\it Rajkov completion $\widehat G$\/} of
a topological group $G$ is the completion of $G$
with respect to the upper uniformity $\sL\vee\sR$.
For every topological group $G$ the space $\widehat G$ has a natural
structure of a topological group. The group $\widehat G$
can be defined as a unique (up to an isomorphism) complete group containing
$G$ as a dense subgroup.
A  group  is  {\it
Weil-complete\/}  if  it  is complete with respect to the left uniformity $\sL$
(or, equivalently, with respect to the right uniformity $\sR$).
Every Weil-complete group is complete, but not vice versa.

Unlike the category of Hausdorff spaces, where ``minimal" implies  ``absolutely
closed",  minimal  groups need not be absolutely closed (that is, complete). If
$G$ is a minimal group, then its  Rajkov  completion  $\widehat  G$
also is minimal. On the other
hand, if $G$ is a dense subgroup of a minimal group $H$, then $G$ is minimal if
and only if for every closed normal subgroup $N\ne\{1\}$ of $H$ we have  $G\cap
N\ne\{1\}$ (\cite{Banache, ProdMin, Stephen}; see historical remarks
in \cite[Section~2.1]{DikSurv}).
Thus the study of minimal groups can be reduced to the
study of complete minimal groups: a group $G$ is minimal if and only if
its Rajkov completion $\widehat G$ is minimal, and for every closed normal subgroup
$N\ne\{1\}$ of $\widehat G$ we have $G\cap N\ne\{1\}$. Compact groups are  complete
minimal, and in the Abelian case the converse is also true, according to a deep
theorem of Prodanov and Stoyanov \cite{PS, DPS}:
every complete minimal Abelian group is
compact. In the non-Abelian case, the class of complete minimal groups properly
contains  the  class  of  compact  groups.  There exist non-compact minimal Lie
groups \cite{DierSchwang, RemusSto},
and actually a discrete infinite group can be minimal \cite{Hesse,Olshan}.  It
is natural to ask how big the difference is between the class of compact groups
and the class of complete minimal groups. For example, one can ask if the class
of  complete  minimal  groups  is  closed  under  infinite  products
(this question, to the best of my knowledge, is still open; the answer is positive
for groups with a trivial center \cite{Meg95}), or if the
relations between cardinal  invariants  of  compact  groups  remain  valid  for
complete minimal groups, etc. 

If $G$ is a topological group, we denote by $\sN(G)$ the filter of 
neighbourhoods of the neutral element.
Besides the left, right, and upper uniformities (denoted by $\sL$, $\sR$,
and $\sL\vee\sR$, respectively), every topological group has yet another 
compatible uniformity $\sL\wedge\sR$, the greatest lower bound of $\sL$ and
$\sR$. (Note that in general the greatest
lower bound of two compatible uniformities on a topological  space
need not be compatible with the topology.)
If $U\in \sN(G)$,  the  cover  $\{UxU:  x\in  G\}$  is
$\sL\wedge\sR$-uniform,  and  every  $\sL\wedge\sR$-uniform  cover of $G$ has a
refinement of this form. 
The uniformity $\sL\wedge\sR$ is called  the  {\it lower uniformity\/} 
in \cite{RD}; we shall call it the {\it Roelcke uniformity\/},
in honour of Walter Roelcke who was
the first to introduce and investigate this notion.

A  uniform  space  $X$  is  {\it
precompact\/} if its completion is  compact  or,  equivalently,  if  for  every
entourage  $U$  the space $X$ can be covered by finitely many $U$-small sets. A
topological group $G$ is {\it precompact\/} if one of the following equivalent
conditions holds\footnote
{Answering a question of
Walter Roelcke, I proved that these conditions are also equivalent
to: (5) for every $U\in\sN(G)$ there 
exists a finite set $F\sbs G$ such that $FUF=G$.
This was later rediscovered by S. Solecki and other authors.
A short proof can be found
in \cite[Proposition 4.3]{BouzTroa}.}: 
(1) $(G,\sL)$ is  precompact; (2) $(G,\sR)$ is precompact;
(3) $(G,\sL\vee\sR)$ is precompact; (4) for  every $U\in\sN(G)$
there  exists  a  finite set $F\sbs G$ such that $FU=UF=G$.
Every  Tikhonov  space  is  a  subspace  of  a  compact  space,  but  not every
topological group is a subgroup of a compact group: the  subgroups  of  compact
groups   are   precisely  precompact  groups.  
Let us say that  a
topological  group  $G$ is  {\it  Roelcke-precompact\/}  
if  it  is precompact with
respect to the Roelcke uniformity $\sL\wedge\sR$. Thus $G$  is
Roelcke-precompact iff for every $U\in\sN(G)$ there exists a finite
set $F\sbs G$ such that $UFU=G$.    
The {\it Roelcke completion\/} of a topological group 
$G$ is the completion
of $G$ with respect to the Roelcke uniformity $\sL\wedge\sR$.
If $G$ is Roelcke-precompact, 
the Roelcke completion $R(G)$ of $G$ will be also called
the {\it Roelcke  compactification\/}.

Precompact groups are Roelcke-precompact, but not  vice  versa
\cite{RD}.  For  example,  the  unitary  group  of  a  Hilbert  space  or  the group
$\operatorname{Sym}(E)$ of all permutations of a discrete set $E$,  
both considered with
the  pointwise convergence topology, are Roelcke-precompact but not precompact.
While the left, right and upper  uniformities  of  a  subgroup  of  a
topological  group  are induced by the corresponding uniformities of the group,
this  is  not  so  for  the  Roelcke  uniformity,   and   a   subgroup of a
Roelcke-precompact group need not be Roelcke-precompact. This justifies
Question~\ref{q:roel}.

The aim of the present paper is to provide a complete answer to
Questions~\ref{q:PesArh} and~\ref{q:roel}.
%
%
Let  us  say  that  a  group  is $G$ {\it topologically simple\/} if $G$ has no
closed normal subgroups besides $G$ and $\{1\}$.

\begin{mainth}
\label{th:main}\label{th:1.1}
Every topological group $G$  is  isomorphic  to  a
subgroup of a complete minimal group which is Roelcke-precompact, topologically
simple and has the same weight as $G$.
\end{mainth}

``Isomorphic"  means  here  ``isomorphic  as  a  topological  group".  The {\it
weight\/} of a topological space  $X$  is  the  cardinal  $w(X)=\min\{|\sB|:\sB
\text{  is  a  base  for  } X\}$. A group $G$ is {\it totally minimal\/} 
\cite{DP} if all
Hausdorff quotient groups
of $G$ are  minimal.  Since  minimal  topologically  simple
groups  are  totally  minimal,  we  could  write  ``totally minimal" instead of
``minimal" in our Main Theorem.

%
%
%
%

Let $Q=[0,1]^\o$ be the Hilbert cube, and let $\Homeo(Q)$ be the topological group
of all self-homeomorphisms of $Q$. 
The group $H=\Homeo(Q)$ is universal \cite{UspFA}, \cite[Theorem 2.2.6]{P1},
in the sense
that every topological group $G$ with a countable base is isomorphic to a topological
subgroup of $H$. Therefore,
for groups with a countable base a natural way to prove Theorem~\ref{th:main}
would be to prove that the group $\Homeo(Q)$, which is known to be simple, is 
Roelcke-precompact and minimal. I do not know if $\Homeo(Q)$ indeed has these properties:

\begin{problem}
\label{pro:minim}
Is the group $\Homeo(Q)$ Roelcke-precompact or minimal?
\end{problem}

There is another universal topological group with a countable base,
namely the group $\Iso(\U)$ of all self-isometries of the Urysohn universal
metric space $\U$ \cite{UspUry}, \cite[Theorem 2.3.1]{P1}, 
\cite[Theorem 6.1]{UspComp}. The group $\Iso(\U)$ is not Roelcke-precompact
\cite[p.\,344]{UspComp};
I do not know whether it is minimal or not. 

\begin{problem}
\label{pro:U}
Is the group $\Iso(\U)$ minimal?
\end{problem}

We consider the ``bounded version"
$\U_1$ of the space $\U$ and show that the group $\Iso(\U_1)$ is 
Roelcke-precompact, topologically simple and minimal. This proves
Theorem~\ref{th:main} for groups with a countable base. For groups
of uncountable weight the argument is similar, but we must consider
non-separable analogues of the space $\U_1$.

Recall some definitions.
A bijection between metric spaces is an {\it isometry\/} if it is
distance-preserving.
For a metric space $M$ we denote by $\Is (M)$ the topological
group of all isometries
of $M$ onto itself, equipped with the topology of pointwise convergence
(which coincides in this case with the compact-open topology).
Let $d$ be the metric on $M$. The sets of the form
$U_{F,\e}=\{g\in\Is(M): d(g(x),x)<\e\text{ for all $x\in F$}\,\}$,
where $F$ is a finite subset of $M$ and $\e>0$, constitute a base at the unity
of $\Is(M)$.

A metric space $M$ is {\it $\o$-homogeneous\/} if
every isometry $f:A\to B$ between finite subsets
$A, B$ of $M$ can be extended to an isometry of $M$ onto itself.
The Urysohn universal metric space $\U$ is the unique (up to an isometry)
complete separable metric space which has either of the following two properties:
(1) $\U$ is $\o$-homogeneous and contains an isometric copy of every separable
metric space; (2) $\U$ is {\em finitely injective}: 
if $L$ is a finite metric space,  $K\sbs  L$  and
$f:K\to  \U$ is an isometric embedding, then $f$ can be extended to an isometric
embedding of $L$ into $\U$. For the equivalence of the conditions (1) and (2),
see Proposition~\ref{p:newequiv} below (we consider there the bounded version
$\U_1$ of $\U$, but the proof for $\U$ is the same). Actually $\U$ is {\em compactly
injective} as well: in the definition of finite injectivity, one can replace finite
metric spaces $K\sbs L$ by arbitrary compact metric spaces \cite{Huhu},
see also \cite[Lemma 5.1.19 and Proposition 5.1.20]{Pbook2}.

We now introduce the bounded version of the space $\U$.
The {\em diameter} of a metric space $(M,d)$ is the number
$\sup\{d(x,y):x,y\in M\}$.
Let us say that a metric space $M$ is {\it Urysohn} if its diameter
is equal to $1$ and it is 
{\em finitely injective with respect to spaces of diameter $\le1$},
that is, the following holds:
if $L$ is a finite metric space of diameter  $\le1$,  $K\sbs  L$  and
$f:K\to  M$ is an isometric embedding, then $f$ can be extended to an isometric
embedding of $L$ into $M$. 
It suffices if this property holds for
$L=K\cup\{p\}$. Thus a metric space $M$ of diameter $1$
is  Urysohn  iff  for  any  finite
sequence  $x_1, \dots, x_n$ of points of $M$ and any sequence $a_1, \dots, a_n$
of positive numbers $\le1$ such that $|a_i-a_j|\le  d(x_i,  x_j)  \le  a_i+a_j\
(i,j=1,\dots,n)$ there exists $y\in M$ such that 
$d(y,x_i)=a_i \ (i=1,\dots,n)$. Using the notion of a Kat\v etov function that
will be introduced later in Section~\ref{s:katconstr}, we can reformulate this 
condition as follows: for every finite $X\sbs M$ and every Kat\v etov function
$f:X\to [0,1]$ there exists $y\in M$ such that $d(x,y)=f(x)$ for every $x\in X$.

\smallskip

{\it Remark.} A notation like Urysohn$_{\le1}$ might have been more appropriate
for what we have called Urysohn (note that the unbounded space $\U$ is not
Urysohn according to our definition!). However, we shall use the shorter term,
in hope that no confusion will arise. Let us again bring to the reader's 
attention that {\em all Urysohn spaces have diameter $1$}.

\begin{prop}
\label{p:newequiv}
Let $M$ be a metric space of diameter $1$.
\begin{enumerate}
\item if $M$ is Urysohn, then $M$ contains an isometric copy of every
countable metric space of diameter $\le1$. If $M$ moreover is
complete, then it contains an isometric copy of every separable
metric space of diameter $\le1$;
\item if $M$ is $\o$-homogeneous and contains an isometric copy
of every finite metric space of diameter $\le1$, then $M$ is Urysohn;
\item if $M_1$ and $M_2$ are complete separable Urysohn spaces, then
every isometry between finite subsets $A\sbs M_1$ and $B\sbs M_2$
extends to an isometry between $M_1$ and $M_2$;
\item a complete separable metric space of diameter 1 is Urysohn
if and only if it is $\o$-homogeneous and contains an isometric
copy of every finite metric space of diameter $\le1$;
\item there exists a unique (up to an isometry) complete separable
Urysohn space $\U_1$. The space $\U_1$ is the unique
complete separable metric space of diameter $\le1$ which is
$\o$-homogeneous and contains
an isometric copy of every separable
metric space of diameter $\le1$;
\item there exists a non-complete separable $\o$-homogeneous Urysohn
space which contains an isometric copy of every separable metric 
space of diameter $\le1$.
\end{enumerate}
\end{prop}

This is essentially due to Urysohn \cite{Ur}. The last item was
added by Kat\v etov \cite{Kat}, 
who answered a question of Urysohn that had remained
open for more than 60 years.

\begin{proof}
(1) is obvious (use induction). To prove (2), suppose that $K\sbs L$
are finite metric spaces, $\diam(L)\le 1$, 
and let $f:K\to M$ be an isometric embedding. Pick an isometric embedding
$g:L\to M$, and use $\o$-homogeneity of $M$ to find an isometry $h$ of $M$
such that $h$ extends the isometry $f(g_{|K})\obr:g(K)\to f(K)$. Then $hg:L\to M$
is an isometric embedding that extends $f$. For (3), enumerate dense
countable subsets in $M_1$ and $M_2$ and use the ``back-and-forth" (or
``shuttle") method to extend the given isometry between $A$ and $B$
to an isometry between dense subsets of $M_1$ and $M_2$. Then use 
completeness to obtain an isometry between $M_1$ and $M_2$.
Applying (3) in the case when $M_1=M_2$, we see 
that every complete separable Urysohn space is
$\o$-homogeneous. Thus (4) and uniqueness in (5) follow from (1)--(3).
The existence of $\U_1$ is a special case of Theorem~\ref{th:UrExt} 
that we shall prove later; the idea of the proof is due to Kat\v etov.
The existence of a non-complete Urysohn space easily follows from
Kat\v etov's methods presented in this paper; we refer the reader to
\cite{Kat} for details.
\end{proof}

For the history of invention of the universal Urysohn space, 
see \cite{PSAlex}, \cite{Tikh}, \cite{Husek}. According to
P.S. Alexandrov \cite{PSAlex}, P.S. Urysohn was thinking about
the universal space in the very last days of his life, and,
after finishing another project on 14~August 1924, was going
to work on two further papers: on metrization of normal spaces
with a countable base and on the universal space. He wrote
just the first page of the first of these papers. It was dated
17~August 1924, the day of his death.

For more on the Urysohn space, see 
\cite{GaoKechr, Mell1, Mell2, P1, P3, Pbook1, Pbook2, UspComp,
UryHilb}, and papers in this volume.
We mention the striking result of Vershik: for a generic
point $d$ of the Polish space of metrics on a countable set $X$
the completion of $(X,d)$ is isometric to the Urysohn space $\U$
\cite{Versh98, Versh02, Versh04}. Similarly, for a generic
shift-invariant metric $d$ on $\Z$ (= the group of integers) 
the completion of the metric group $(\Z,d)$ is isometric to $\U$
\cite{CamVersh}.

The proof of Theorem~\ref{th:main} consists of two parts. 
We first prove that every
topological group can be embedded in the group $\Is(M)$
of isometries of a complete $\o$-homogeneous Urysohn space $M$, 
and then prove that such groups of isometries
are minimal, Roelcke-precompact and topologically simple. 
\begin{thm}
\label{th:1.4}
For every topological group $G$ there exists a
complete $\o$-homogeneous
Urysohn metric space $M$ of the same weight as $G$ such that $G$
is isomorphic to a subgroup of $\Is(M)$.
\end{thm}

\begin{thm}
\label{th:1.5}
If $M$ is a complete $\o$-homogeneous
Urysohn metric space, then the group $\Is(M)$ is complete, Roelcke-precompact,
minimal and topologically simple. The weight of $\Is (M)$ is equal to the
weight of $M$.
\end{thm}

Theorem~\ref{th:1.1} follows from Theorems~\ref{th:1.4} and~\ref{th:1.5}.

The proof of Theorem~\ref{th:1.4} depends on Kat\v etov's construction
that leads to a canonical embedding of any metric space $M$ into a finitely
injective space. In the non-separable case this construction must
be complemented by a construction of a natural embedding of a metric
space into an $\o$-homogeneous space. We use Graev metrics on free groups
for this.

The proof of Theorem~\ref{th:1.5} is based on
the study of the Roelcke compactifications of groups of isometries.
The Roelcke compactifications of some topological groups of importance
admit an explicit description and are equipped with additional structures.
For example, for the unitary group $U_s(H)$, where $H$ is a Hilbert space
and the subscript $s$ indicates the strong operator topology (= the topology
inherited from the Tikhonov product $H^H$), the Roelcke compactification
can be identified with the unit ball $\Theta$ in the algebra of bounded
linear operators on $H$ \cite{Orsat}. The ball $\Theta$ is equipped
with the weak operator topology. This is the topology inherited from 
$H^H$, where each factor $H$ carries the weak topology. 
Another case when the Roelcke compactification can be explicitly described
is the following. Let $K$ be a zero-dimensional compact space such that
all non-empty clopen subspaces of $K$ are homeomorphic to $K$. For example,
$K$ might be the Cantor set. Let $G=\Homeo(K)$, equipped with the 
compact-open topology. Then 
$R(G)$ is the 
set of all closed relations $R$ on $K$ (= closed subsets of $K^2$) 
such that the domain and the range of $R$ is equal to $K$ \cite{U4}. 
Yet another
example of a topological group $G$ for which $R(G)$ 
is known is the group $G=\Homeo_+[0,1]$ of all orientation-preserving
self-homeomorphisms of the closed interval $I=[0,1]$. In that case
$R(G)$ can be identified with the closure of the
set of graphs of elements of $G$ in the space of closed
subsets of the square $I^2$, see the picture in 
\cite[Example 2.5.4]{P1}.

The proof of Theorem~\ref{th:1.5} leans on the study of
the Roelcke compactification $R(G)$ for $G=\Is(M)$,
where $M$ is a complete $\o$-homogeneous Urysohn metric space.
In this case $R(G)$ can be identified with the space of metric
spaces covered by two isometric copies of $M$, see 
Sections~\ref{s:semi} and~\ref{s:roelcomp} below. 
Equivalently, $\Th=R(G)$ can be identified with 
a certain subset of $I^{M\times M}$ that we now are going
to specify.

A {\it semigroup\/} is a set with an associative binary operation.
Let $S$ be a semigroup with the multiplication $(x,y)\mapsto xy$.
An element $x\in S$ is an {\it idempotent\/} if $x^2=x$.
We say that a self-map
$x\mapsto x^*$ of $S$ is an {\it involution\/} if $x^{**}=x$ and
$(xy)^*=y^*x^*$ for all $x,y\in S$.  
An element $x\in S$ is {\it symmetrical\/} if $x^*=x$, and a subset $A\sbs S$
is {\it symmetrical\/} if $A^*=A$.
An {\it ordered semigroup\/} is a semigroup with a partial order $\le$ such
that the conditions $x\le x'$ and $y\le y'$ imply $xy\le x'y'$.

Denote by $I$ the closed unit interval $[0,1]$. Let $\uplus$ be the
associative operation on $I$ defined by $x\uplus y=\min(x+y,\,1)$. 
Let $X$ be a set, and let $S=I^{X\times X}$ be the set of all functions 
$f:X^2\to I$. 
We make $S$ into an ordered semigroup with an involution. 
Define an operation $(f,g)\mapsto f\bull g$ on
$S$ by
$$
f\bull g(x,y)=\inf\{f(x,z)\uplus g(z,y):z\in X\} \quad (x,y\in X).
$$
This operation is associative, since for $f,g,h\in S$ and 
$x,y\in X$ both $(f\bull g)\bull h(x,y)$ $f\bull (g\bull h)(x,y)$ are equal to
$$
\inf\{f(x,z)\uplus g(z,u)\uplus h(u,y):z,u\in X\}.
$$
Define an involution $f\mapsto f^*$ on $S$ by $f^*(x,y)=f(y,x)$. 

Let $(M,d)$ be a complete $\o$-homogeneous Urysohn metric space, 
and let $G=\Is(M)$.
The Roelcke compactification $\Th$ of $G$ 
can be identified with a closed subsemigroup of $I^{M\times M}$
and has a natural structure of an ordered
semigroup with an involution. Namely, $\Th$ can be viewed as
the set of all functions $f\in I^{M\times M}$ 
which are bi-Kat\v etov in the sense of Definition~\ref{def:bikat}.
Such functions can be described in terms of the structure
of an ordered semigroup with an involution 
on $I^{M\times M}$: a function $f\in I^{M\times M}$ is bi-Kat\v etov
if and only if
$$
f\bull d=d\bull f=f, \quad f^*\bull f\ge d, \quad f\bull f^*\ge d.
$$
The metric $d$ is the unity of $\Th$, and the constant 1 is a zero element
of $\Theta$, in the sense that $f\bull 1=1\bull f=1$ for every $f\in \Theta$
(in fact, for every $f\in I^{M\times M}$).

Note that $\Theta$ is a compact space and a semigroup, but it might be
misleading to call it a ``compact semigroup", since the semigroup operation
on $\Theta$ is not (even separately) continuous. However, both the topology
and the algebraic structure on $\Theta$ will play an important role
in our proofs.

The Roelcke compactification $\Theta$ of $G=\Is(M)$ is used 
to prove Theorem~\ref{th:1.5} in the following way.
Let $f:G\to G'$ be a surjective morphism of topological groups. To prove that
$G$ is minimal and topologically simple, we must prove that either $f$ is
a homeomorphism or $|G'|=1$. Extend $f$
to a map $F:\Th\to \Th'$, where $\Th'$ is the Roelcke compactification of $G'$.
Let $S=F\obr(e')$, where $e'$ is the unity of $G'$. Then $S$ is a closed 
subsemigroup of $\Th$ which is invariant under inner automorphisms. To every 
closed subsemigroup of $\Th$ an idempotent can be assigned in a canonical way.
Let $p$ be the idempotent corresponding to $S$. Since $S$ is invariant under
inner automorphisms, so is $p$. We show that certain idempotents in $\Th$ 
are in a one-to-one correspondence with closed subsets of $M$
(Proposition~\ref{p:idemp}). Since there are no non-trivial $G$-invariant
closed subsets of $M$,
it follows that $p$ is trivial: it is either the unity of
$\Th$ or the constant~1. Accordingly, either $f$ is a homeomorphism or
$G'=\{e'\}$.

The same 
method was used in \cite{Orsat} and \cite{U4}
to give alternative proofs of Stoyanov's theorem that the unitary group of
a Hilbert space is totally minimal and of Gamarnik's theorem
that the group of homeomorphisms of the Cantor set is minimal,
see Remarks~2 and~3 in Section~\ref{s:rem}
below.

Under the conditions of Theorem~\ref{th:1.5}, the group $\Is(M)$ has the
{\it fixed point on compacta (f.p.c.)\/} property. This deep result is due
to V.~Pestov \cite{P3, PesLevy, Pbook1, Pbook2}%
\footnote{
The setting considered in these papers and books is not exactly the same as in
Theorem~\ref{th:1.5} (detailed proofs are given either for the separable case
or for unbounded metrics), but, as noted in \cite{P3}, the same argument works
for bounded metrics {\it verbatim}.}.
A topological group $G$ has the f.p.c.
\hskip -1pt
property, or is {\em extremely amenable},
if every compact $G$-space has a $G$-fixed point.
As pointed out by Pestov, his theorem, combined with Theorem~\ref{th:1.4}
of the present paper, implies that every topological group is a subgroup
of an extremely amenable group.

We prove Theorem~\ref{th:1.4}
in Section~\ref{s:proof1.4} and Theorem~\ref{th:1.5} in 
Section~\ref{s:proof1.5}.

Another version of Question~\ref{q:PesArh} is the following
(see \cite[Problem~VI.6]{Arh}, \cite[Problem~519]{OpenPr}):
is every topological group a quotient
of a minimal topological group?
I have earlier announced that the answer is positive.
Moreover, I claimed that every topological group is a group
retract of a minimal topological group. 
In other words, for every topological group $G$ there exist
a minimal topological group $G'\supseteq G$ and
a morphism $r:G'\to G$ such that $r^2=r$
(it follows that $G$ is a quotient of $G'$).
My announcement appears as Theorem~3.3F.2 in \cite{CHR}.
However, my announcement was premature, and my ``proof" contained a gap.
A complete proof has been recently found by M. Megrelishvili \cite{Megr}.
%

Megrelishvili's construction shows that every complete group is a group
retract of a complete minimal group. This result, combined with the fact 
that every topological group is a quotient of a Weil-complete group 
\cite{UspIzv}%
\footnote{It was proved in \cite{UspIzv} that the free topological group 
of any stratifiable space is Weil-complete. Since every topological space 
is the image of a stratifiable space under a quotient (even open) map 
\cite{Jun},
it follows that every topological group is a quotient of a Weil-complete group.},
implies that {\em every topological group is a quotient of a complete minimal group.}
Indeed, given any topological group $G$, represent $G$ as a quotient of a complete
group $G'$, and then, using Megrelishvili's theorem,
represent $G'$ as a group retract (and hence as a quotient)
of a complete minimal group.

\section{Invariant pseudometrics on groups}
\label{s:2}

A pseudometric $d$ on a group $G$ is {\it left-invariant\/}  if
$d(xy, xz)=d(y,z)$ for all $x,y,z\in G$. Right-invariant pseudometrics
are defined similarly. A pseudometric is {\it two-sided invariant\/} if
it is left-invariant and right-invariant. Let $e$ be the unity of $G$.
A non-negative real function $p$ on $G$ is a {\it seminorm\/} if it satisfies
the following conditions: (1) $p(e)=0$; (2) $p(xy)\le p(x)+p(y)$ for all
$x,y\in G$; (3) $p(x\obr)=p(x)$ for all $x\in G$. If $p$ is a seminorm on $G$,
define a left-invariant pseudometric $d$ by $d(x,y)=p(x\obr y)$.
We thus get a one-to-one correspondence between seminorms and left-invariant
pseudometrics.
Given a left-invariant pseudometric $d$, the corresponding seminorm $p$ is
defined by $p(x)=d(x,e)$. A seminorm $p$ is {\it invariant\/} if it is
invariant under inner automorphisms, that is $p(yxy\obr)=p(x)$ for every
$x,y\in G$. Invariant seminorms correspond to two-sided invariant
pseudometrics.

Now let $G$ be a topological group. Then the topology of $G$ is determined
by the collection of all continuous left-invariant pseudometrics
\cite[Theorem~8.2]{HR}. Equivalently,
for every neighbourhood $U$ of unity there exists a continuous seminorm $p$
on $G$ such that the set $\{x\in G: p(x)<1\}$ is contained in $U$.

\begin{thm}
\label{th:GrIsom}
For every topological group $G$ there exists a metric space $M$ such that
$w(G)=w(M)$ and
$G$ is isomorphic (as a topological group) to a subgroup of $\Is(M)$.
\end{thm}

This theorem has been rediscovered many times by various authors,
see historical remarks in \cite{P1, P2}.

\begin{proof}[1st proof]
There exists a family $D=\{d_\a:\a\in A\}$ of continuous left-invariant
pseudometrics on $G$ which determines the topology of $G$ and has the
cardinality $|A|=w(G)$. Replacing, if necessary, each $d\in D$ by
$\inf(d,1)$, we may assume that all pseudometrics in $D$ are bounded by~1.
For every $\a\in A$ let $M_\a$ be the metric space associated with the
pseudometric space $(G,d_\a)$, and let $M=\bigoplus_{\a\in A}M_\a$ be the
disjoint sum of the spaces $M_\a$. There is an obvious metric on $M$
which extends the metric of each $M_\a$:
if two points of $M$ are in distinct pieces $M_\a$ and $M_\beta$, define
the distance between them to be~1.
The left action of $G$ on itself yields for every $\a\in A$
a natural continuous homomorphism $G\to \Is(M_\a)$.
The homomorphism $G\to \prod_{\a\in A}\Is(M_\a)$ thus obtained is
a homeomorphic embedding. It remains to note that the group
$\prod_{\a\in A}\Is(M_\a)$  can be identified with a topological subgroup of
$\Is(M)$.
\end{proof}

\begin{proof}[2nd proof]
Let $B$ be the Banach space of all bounded real functions on $G$ which are
uniformly continuous with respect to the right uniformity. The natural left
action of $G$ on $B$, defined by the formula $gf(h)=f(g\obr h)$
$(g,h\in G,\ f\in B)$, yields an isomorphic embedding of $G$ into $\Is(B)$.
The weight of $B$ may exceed the weight of $G$, but it is easy to find
a $G$-invariant subspace $B'$ of $B$ such that $B'$ determines
the topology of $G$ and $w(B')=w(G)$. Then the natural
homomorphism $G\to \Is(B')$ still is a homeomorphic embedding.
\end{proof}

Let us discuss invariant seminorms on free groups.
For a set $X$ we denote by $S(X)$ the set of all words of the form
$x_1^{\e_1}\dots x_n^{\e_n}$, where $n\ge 0$, $x_i\in X$ and $\e_i=\pm1$,
$1\le i\le n$. In other words, $S(X)$ is the free monoid
\footnote
{A {\em monoid} is a semigroup with a neutral element. We require that monoid
morphisms should preserve the neutral element.}
on the set
$X\cup X\obr$, where $X\obr$ is a disjoint copy of $X$. A word $w\in S(X)$
is {\it irreducible\/} if it does not contain subwords of the form
$x^\e x^{-\e}$.
We consider the {\it free group\/} $F(X)$ on a set $X$ as the set of
all irreducible words in $S(X)$. Every word $w\in S(X)$ represents a uniquely
defined element $w'\in F(X)$ which can be obtained from $w$ by
consecutive deletion of subwords of the form $x^\e x^{-\e}$.
In this situation we say that the words $w$ and $w'$ are {\it equivalent}.
For
$u,v\in S(X)$ we denote by $u|v$ the product of $u$ and $v$ in the semigroup
$S(X)$, that is the word obtained by writing $v$ after $u$
(without cancelations). If $u$ and $v$ are irreducible, we denote by $uv$
their product in the group $F(X)$, that is the irreducible word equivalent to
$u|v$.

Let $(X,d)$ be a metric space. A real function $f$ on $X$ is
{\it non-expanding\/}, or {\em $1$-Lipschitz},
if $|f(x)-f(y)|\le d(x,y)$ for every $x,y\in X$.
Let $k$ be a non-negative non-expanding function on $X$. We shall
describe a two-sided invariant pseudometric $Gr(d,k)$ on the free group
$F(X)$ which is called the {\it Graev pseudometric}
\cite{Graev}. The corresponding invariant
seminorm $p$ is characterized by the following property: $p$ is the greatest
invariant seminorm on $F(X)$ such that $p(x)=k(x)$ and $p(x\obr y)\le d(x,y)$
for every $x,y\in X$. We shall need later
the following explicit construction of the seminorm $p$.

It will be convenient to define the function $p$ on the entire set $S(X)$.
Given a word $w=x_1^{\e_1}\dots x_n^{\e_n}\in S(X)$, we define a
{\it $w$-pairing\/} to be a collection $E$ of pairwise disjoint two-element
subsets of the set $J=\{1,\dots,n\}$ such that: (1) if $\{a,b\}\in E$ and
$\{i,j\}\in E$, where $a<b$ and $i<j$, then the intervals $[a,b]$ and $[i,j]$
are either disjoint or one of them is contained in the other (this means that
the cases $a<i<b<j$ and $i<a<j<b$ are excluded); (2) if $\{i,j\}\in E$,
then $\e_i=-\e_j$. To put it less formally, some pairs of letters of the word
$w$ are connected by arcs, each letter is connected with at most one other
letter, each arc connects a pair of letters of the form $x$ and $y\obr$
($x,y\in X$), and the arcs do not intersect each other.
Given a $w$-pairing $E$, define the {\it Graev sum\/}
$s_E=s_E(w)$ by
$$
s_E =\sum \{d(x_i,x_j): \{i,j\}\in E,\ i<j\} +
 \sum \{k(x_i):i\in J\setminus \cup E\},
$$
and let $p(w)$ be the minimum of the numbers $s_E$, taken over the finite
set of all $w$-pairings $E$.

We claim that $p(w)=p(w')$ if the words $w,w'\in S(X)$ are equivalent. It
suffices to consider the case when $w=u|v$ and $w'=u|x^\e x^{-\e}|v$.
We show that for every $w'$-pairing $E'$ there exists a $w$-pairing $E$
such that $s_E\le s_{E'}$, and vice versa. In one direction this is obvious:
given a $w$-pairing $E$, which we consider as a system of arcs connecting
the letters of the word $w$, add one more arc which connects the letters
$x^\e$ and $x^{-\e}$ of the word $w'$. We get a $w'$-pairing $E'$ for which
$s_E=s_{E'}$. Conversely, let a $w'$-pairing $E'$ be given. We must construct
a $w$-pairing $E$ for which $s_E\le s_{E'}$. As above, we consider $E'$ as
a system of arcs. The word $w$ is obtained from $w'$ by deleting the subword
$x^\e x^{-\e}$. To get $E$, we replace the arcs which go from the letters
$x^\e$ and $x^{-\e}$ and leave the other arcs unchanged.
Consider the following cases.

Case 1. There is an arc in $E'$ connecting the letters $x^\e$ and $x^{-\e}$.
Then just delete this arc to get $E$. We have $s_E=s_{E'}$.

Case 2. The letters $x^\e$ and $x^{-\e}$ are connected in $E'$, but not with
each other. Let $x^\e$ be connected with $y^{-\e}$ and $x^{-\e}$ be connected
with $z^\e$. Replace these two connections by one connection between $y^{-\e}$
and $z^\e$. The sums $s_E$ and $s_{E'}$ differ by the terms $d(y,z)$ and
$d(y,x)+d(x,z)$, hence the triangle inequality implies that $s_E\le s_{E'}$.

Case 3. One of the letters $x^\e$ and $x^{-\e}$, say $x^\e$, is connected
in $E'$ and the other is unpaired. Let $x^\e$ be connected with $y^{-\e}$.
Delete this connection and leave the letter $y^{-\e}$ unpaired in $E$.
The sums $s_E$ and $s_{E'}$ differ by the terms $k(y)$ and $d(x,y)+k(x)$.
Since the function $k$ is non-expanding, we have $k(y)\le d(x,y)+k(x)$
and hence $s_E\le s_{E'}$.

Case 4. Both $x^\e$ and $x^{-\e}$ are unpaired in $E'$. Then all arcs are
left without change. The sum $s_E$ is obtained from $s_{E'}$ by omitting
the non-negative term $2k(x)$, hence $s_E\le s_{E'}$.

We have thus proved the claim that
$p(w)=p(w')$ for equivalent words $w,w'\in S(X)$.
It follows that the restriction of $p$ to $F(X)$ is indeed a seminorm:
if $u,v\in F(X)$, then $p(uv)=p(u|v)\le p(u)+p(v)$. It is easy to see that
$p(u)=p(u\obr)$ for every $u\in F(X)$. We show that $p$ is invariant under
inner automorphisms. If $u\in S(X)$, $x\in X$, $\e=\pm1$ and
$w=x^\e|u|x^{-\e}$, then
$p(w)\le p(u)$, since every $u$-pairing can be extended 
in an obvious way to a $w$-pairing
with the same Graev sum.
It follows that for every $u,v\in F(X)$
we have $p(uvu\obr)=p(u|v|u\obr)\le p(v)$, and by symmetry 
of the relation of being conjugate in $F(X)$ also the opposite
inequality holds. Thus $p(uvu\obr)=p(v)$, which means that the seminorm $p$
is invariant.

Let $Y$ be a pseudometric space, and let $\Is(Y)$ be the group of all
distance-preserving permutations of $Y$, equipped with 
the topology of pointwise 
convergence. Then $\Is(Y)$ is a topological group, not necessarily Hausdorff.
For later use we note here the following:

\begin{lemma}
\label{l:Gr-metr-isom}
Let $(X,d)$ be a metric space, and let $k$ be a non-expanding function on $X$.
Let $D=Gr(d,k)$ be the Graev pseudometric on the free group $G=F(X)$. Let
$H_1\sbs\Is(X)$
be the topological group of all isometries of $X$ which preserve the function
$k$, and let $H_2\sbs\Is(G)$ be the topological group (not necessarily
Hausdorff) of all automorphisms of $G$ which preserve the pseudometric
$D$. Then the natural homomorphism $\f\mapsto\f^*$
from $H_1$ to $H_2$ is continuous.
\end{lemma}

\begin{proof}
It suffices to show that for every $w\in G$ the map $\f\mapsto\f^*(w)$
from $H_1$ to $(G,D)$ is continuous at the unity.
If $w=x_1^{\e_1}\dots x_n^{\e_n}$, then
$\f^*(w)=\f(x_1)^{\e_1}\dots \f(x_n)^{\e_n}$, and
we have $D(\f^*(w), w)\le \sum_{i=1}^n d(\f(x_i),x_i)$. Let $\e>0$ be given.
If $\f\in H_1$ is close enough to the unity, then  $d(\f(x_i),x_i)<\e/n$,
$1\le i\le n$, and therefore $D(\f^*(w), w)<\e$.
\end{proof}

\section
{Kat\v etov's construction of Urysohn extensions}
\label{s:katconstr}

\begin{defin}
\label{def:g-embed}
Let $M$ be a subspace of a metric space $L$. We say that $M$ is 
{\it $g$-embedded\/} in $L$ if there exists a continuous homomorphism
$e:\Is(M)\to\Is(L)$ such that for every $\f\in\Is(M)$
the isometry $e(\f):L\to L$ is an extension of $\f$.
\end{defin}

Let $M$ be a $g$-embedded subspace of a metric space $L$.
A homomorphism $e:\Is(M)\to\Is(L)$ satisfying the condition of
Definition~\ref{def:g-embed}
is a homeomorphic embedding, since the inverse map
$e(\f)\mapsto \f=e(\f)|M$
is continuous. It follows that $\Is(M)$ is isomorphic to a topological subgroup
of $\Is(L)$.

In this section we prove the following theorem:

\begin{thm}
\label{th:UrExt}
Let $M$ be a metric space of diameter $\le1$. There exists a complete
Urysohn metric space $L$ containing $M$ as a subspace such that $w(L)=w(M)$
and $M$ is $g$-embedded in $L$.
\end{thm}

It follows that for every topological group $G$ there exists a
complete Urysohn metric space $M$ of the same weight as $G$ such that $G$
is isomorphic to a subgroup of $\Is(M)$. This is weaker than
Theorem~\ref{th:1.4},
since in the non-separable case
the metric space $M$ need not be $\o$-homogeneous.
In the next section we shall prove that every metric space $M$ can be
$g$-embedded into an $\o$-homogeneous metric space $L$. Using this fact,
we show that the Urysohn space $L$ in Theorem~\ref{th:UrExt} can 
be additionally assumed $\o$-homogeneous (Theorem~\ref{th:pred1.4}).
This yields Theorem~\ref{th:1.4}, see Section~\ref{s:proof1.4}.

The arguments of \cite{UspUry, UspComp} show that Theorem~\ref{th:UrExt}
essentially follows from Kat\v etov's construction of
Urysohn extensions \cite{Kat}.
For the reader's convenience we give a detailed proof.

Let $(X,d)$ be a 
metric space of diameter $\le1$. We say that a function
$f:X\to [0,1]$ is {\it Kat\v etov\/} if $|f(x)-f(y)|\le d(x,y) \le f(x)+f(y)$
for all $x,y\in X$. A function $f$ is Kat\v etov if and only if there exists
a metric space $X'=X\cup \{p\}$ of diameter $\le1$
containing $X$ as a subspace such that 
for every $x\in X$ $f(x)$ is equal to the distance between $x$ and $p$.
Let $E(X)$ be the set of all Kat\v etov functions on $X$, equipped with the
sup-metric $d_X^E$ defined by $d_X^E(f,g)=\sup\{|f(x)-g(x)|:x\in X\}$. If 
$Y$ is a non-empty subset of $X$ and $f\in E(Y)$, 
define $g=\k_Y(f)\in E(X)$ by
$$
g(x)=\inf(\{d(x,y)+f(y):y\in Y\}\cup\{1\})
=\inf\{d(x,y)\uplus f(y):y\in Y\}
$$
for every $x\in X$.
It is easy to check that $g$ is indeed a Kat\v etov function on $X$ and that
$g$ extends $f$; one can define $g$ as 
the largest 1-Lipschitz function $X\to[0,1]$
that extends $f$.
The map $\k_Y:E(Y)\to E(X)$ is
an isometric embedding.
Let
$$
E(X,\o)=
\bigcup\{\k_Y(E(Y)): Y\sbs X,\  Y\text{ is finite and non-empty}\,\}
\sbs E(X).
$$
For every $x\in X$ let $h_x\in E(X)$ be the function on $X$
defined by $h_x(y)=d(x,y)$. Note that $h_x=\k_{\{x\}}(0)$ and hence
$h_x\in E(X,\o)$.
The map $x\mapsto h_x$ is an isometric embedding of $X$ into $E(X,\o)$.
Thus we can identify $X$ with a subspace of $E(X,\o)$. 

\begin{lemma}
\label{l:dist}
If $x\in X$ and $f\in E(X)$, then $d_X^E(f,h_x)=f(x)$.
\end{lemma}

\begin{proof}
Since $f$ is a Kat\v etov function, for every $y\in Y$ we have
$f(y)-d(x,y)\le f(x)$ and $d(x,y)-f(y)\le f(x)$. Hence
$d_X^E(f,h_x)=\sup\{|f(y)-d(x,y)|:y\in X\} \le f(x)$, and at $y=x$
the equality is attained.
\end{proof}

\begin{lemma}
\label{l:yp}
Let $Z=Y\cup\{p\}$ be a finite metric space of diameter $\le1$. Every
isometric embedding $j:Y\to X$ extends to an isometric embedding
of $Z$ into $E(X,\o)$.
\end{lemma}

\begin{proof}
We may assume that $Y$ is a subspace of $X$
and that $j(y)=y$ for every $y\in Y$. Let $f\in E(Y)$ be the
Kat\v etov function defined by $f(y)=\nu(y, p)$ for every $y\in Y$,
where $\nu$ is the metric on $Z$. Let $g=\k_Y(f)\in E(X,\o)$. We claim that
the extension of $j$ over $Z$ which maps $p$ to $g$ is an isometric embedding.
It suffices to check that $d_X^E(h_y, g)=\nu(y,p)$ for every
$y\in Y$. Fix $y\in Y$. Let $h_y^*\in E(Y)$ be the restriction of 
$h_y$ to $Y$. According to Lemma~\ref{l:dist} we have $d_Y^E(h_y^*, f)=f(y)$.
Since $h_y=\k_Y(h_y^*)$, $g=\k_Y(f)$ and the map $\k_Y:E(Y)\to E(X)$
is distance-preserving, it follows that $d_X^E(h_y, g)=d_Y^E(h_y^*, f)
=f(y)=\nu(y, p)$, as claimed.
\end{proof}

\begin{lemma}
\label{l:urg}
Any metric space $X$ of diameter $\le1$ is $g$-embedded in $E(X,\o)$.
\end{lemma}

\begin{proof} 
It is clear that every isometry $\f:Y\to Z$ between any two metric spaces
can be extended to an isometry $\f^*:E(Y,\o)\to E(Z,\o)$. Such an extension
is unique, since every point in $E(Y,\o)$ (or, more generally, in $E(Y)$)
is uniquely determined by its distances from the points of $Y$
(Lemma~\ref{l:dist}), and similarly
for $Z$. In particular,
every isometry $\f\in \Is(X)$ uniquely extends to an isometry
$\f^*\in \Is(E(X,\o))$. The map $\f\mapsto\f^*$ is a homomorphism
of groups. We show that this homomorphism is continuous. Fix $f\in E(X,\o)$
and $\e>0$. Pick a finite subset $Y$ of $X$ and $g\in E(Y)$ so that
$f=\k_Y(g)$. Let $U$ be the set of all $\f\in \Is(X)$ such that
$d(\f(y),y)<\e$ for every $y\in Y$. Then $U$ is a neighbourhood of unity
in $\Is(X)$. It suffices to show that $d_X^E(\f^*(f), f)<\e$ for every
$\f\in U$. Fix $\f\in U$. Let $g_\f=g\circ \f\obr\in E(\f(Y))$.
Then $\f^*(f)=\k_{\f(Y)}(g_\f)$. Thus for every $x\in X$ we have
$$
\f^*(f)(x)=\inf\{d(x,z)\uplus g_\f(z):z\in \f(Y)\}
=\inf\{d(x,\f(y))\uplus g(y):y\in Y\}.
$$
Since
$$
f(x)=\inf\{d(x,y)\uplus g(y):y\in Y\},
$$
it follows that
$$
|\f^*(f)(x)-f(x)|\le \sup\{|d(x,\f(y))-d(x,y)|:y\in Y\}\le
\max\{d(y,\f(y)):y\in Y\}<\e,
$$
whence $d_X^E(\f^*(f), f)<\e$.
\end{proof}

Let $\a$ be an ordinal, and let $\sM=\{M_\beta:\beta<\a\}$ be a family
of metric spaces such that $M_\beta$ is a subspace of $M_\g$ for all
$\beta<\g<\a$. We say that 
the family $\sM$ is {\it continuous\/} if $M_\beta=\bigcup_{\g<\beta}M_\g$
for every limit ordinal $\beta<\a$, $\beta>0$.

\begin{prop}
\label{p:gchain}
Let $\{M_\beta:\beta\le\a\}$ 
be an increasing continuous chain of metric spaces.
If $M_\beta$ is $g$-embedded in $M_{\beta+1}$ for every $\beta<\a$, then 
$M_0$ is $g$-embedded in $M_\a$.
\end{prop}

\begin{proof}
For every $\beta<\a$ pick a continuous homomorphism
$e_\beta:\Is(M_\beta)\to \Is(M_{\beta+1})$ such that $e_\beta(\f)$ extends
$\f$ for every $\f\in \Is(M_\beta)$. 
By transfinite recursion on $\beta\le\a$
define a homomorphism $\La_\beta:\Is(M_0)\to \Is(M_\beta)$ such 
that $\La_\beta(\f)$ extends $\La_\g(\f)$ for every $\f\in \Is(M_0)$ and
$\g<\beta\le\a$. 
Let $\La_0$ be the identity map of $\Is(M_0)$. If
$\beta=\g+1$, put $\La_\beta=e_\g \La_\g$. If $\beta$ is a limit ordinal,
let $\La_\beta(\f)$ be the isometry of $M_\beta$ such that
for every $\g<\beta$ its restriction
to $M_\g$ is equal to $\La_\g(\f)$. We prove by
induction on $\beta$ that each homomorphism $\La_\beta$ is continuous.
This is obvious for non-limit ordinals. Assume that $\beta$ is limit.
To prove that $\La_\beta:\Is(M_0)\to \Is(M_\beta)$ is continuous,
it suffices to show that for every $x\in M_\beta$ the map 
$\f\mapsto \La_\beta(\f)(x)$ from $\Is(M_0)$ to $M_\beta$ is continuous.
Fix $x\in M_\beta$. Pick $\g<\beta$ so that $x\in M_\g$. Then 
$\La_\beta(\f)(x)=\La_\g(\f)(x)$ for every $\f\in \Is(M_0)$.
The map $\La_\g$ is continuous by the assumption of induction, hence
the map $\f\mapsto \La_\beta(\f)(x)=\La_\g(\f)(x)$ also is continuous.
Thus $\La_\a:\Is(M_0)\to \Is(M_\a)$ is a continuous homomorphism such that
$\La_\a(\f)$ extends $\f$ for every $\f\in \Is(M_0)$.
This means that $M_0$ is $g$-embedded in $M_\a$.
\end{proof}     

Put $X_0=X$, $X_{n+1}=E(X_n, \o)$. We consider each $X_n$ as a subspace of
$X_{n+1}$, so we get an increasing sequence $X_0\sbs X_1\sbs\dots$ of
metric spaces. Let $X_\o=\bigcup\{X_n:n\in\o\}$. 

\begin{prop}
\label{p:propur}
The space $X_\o$ is Urysohn, and $X$ is $g$-embedded in $X_\o$.
\end{prop}

\begin{proof}
Let $Z=Y\cup\{p\}$ be a finite metric space of diameter $\le1$, and let 
$j:Y\to X_\o$ be an isometric embedding. Pick $n\in \o$ so that
$j(Y)\sbs X_n$. In virtue of Lemma~\ref{l:yp}, there exists an isometric
embedding of $Z$ into $X_{n+1}\sbs X_\o$ which extends $j$. This means that
$X_\o$ is Urysohn. The second assertion of the proposition follows from
Lemma~\ref{l:urg} and Proposition~\ref{p:gchain}.
\end{proof}

\begin{prop}[\cite{Kat}]
\label{p:Katweight}
The weight of $X_\o$ is equal to the weight of $X$.
\end{prop}

\begin{proof}
It suffices to show that for every metric space $X$
the weight of $E(X,\o)$ is equal to the weight of
$X$. Let $w(X)=\tau$, and let $A$ be a dense subset of $X$ of cardinality
$\tau$. Let $\g=\{\k_Y(E(Y)): Y\sbs A,\  Y\text{ is finite}\}$. Then $\g$
is a family of separable subspaces of $E(X,\o)$, $|\g|=\tau$ and $\cup\g$
is dense in $E(X,\o)$
(see the proof of Lemma~1.8 in \cite{Kat}). Hence $E(X,\o)$ has a dense subspace
of cardinality $\tau$.
\end{proof}

\begin{prop}
\label{p:compl}
Every metric space is $g$-embedded in its completion.
\end{prop}

\begin{proof}
Let $M$ be a metric space, $\cl M$ be its completion. Every isometry
$\f\in \Is(M)$ uniquely extends to an isometry $\f^*\in\Is(\cl M)$.
We show that the homomorphism $\f\mapsto\f^*$ is continuous. Let 
$d$ be the metric on $\cl M$. Fix $x\in \cl M$ and $\e>0$. Pick 
$y\in M$ so that $d(x,y)<\e$. Let $U=\{\f\in \Is(M):d(\f(y),y)<\e\}$.
Then $U$ is a neighbourhood of unity in $\Is(M)$. If $\f\in U$, then
$d(\f^*(x),x)\le d(\f^*(x),\f^*(y))+d(\f^*(y),y)+d(y,x) <3\e$. This 
implies the continuity of the homomorphism $\f\mapsto\f^*$.
\end{proof}

\begin{prop}[{\cite{Ur}}, {\cite[Lemma 3.4.10]{Pbook1}},
{\cite[Lemma 5.1.17]{Pbook2}}, {\cite[Section 3.11$\frac{2}{3}_+$]{Gromov}}]
\label{p:ThUr}
The completion of any Urysohn metric space is Urysohn.
\end{prop}

\begin{proof}
Let $(M,d)$ be a complete metric space containing a dense Urysohn subspace
$A$. We must prove that $M$ is Urysohn.

Let $Y$ be a finite subset of $M$, and let $f\in E(Y)$ be a Kat\v etov
function. It suffices to prove that there exists $z\in M$ such that $d(y,z)
=f(y)$ for every $y\in Y$. Pick a sequence $\{a_n:n\in \o\}\sbs A$ such that:
\begin{itemize}
\item[] if $A_n=\{a_k:k\le n\}$ and $r_n=d(a_{n+1}, A_n)$, $n=0,1,\dots$,
then the series $\sum r_n$ converges;
\item[] every $y\in Y$ is a cluster point of the sequence $\{a_n:n\in \o\}$.
\end{itemize}
To construct such a sequence, enumerate $Y$ as $Y=\{y_1, \dots, y_s\}$, and
for every $k$ and $j$ ($k\in \o$, $1\le j\le s$) pick a point $x_k^j\in A$
such that $d(x_k^j, y_j)<2^{-k}$. Then $d(x_{k+1}^j, x_k^j)<2^{1-k}$ for every
$k$ and $j$, and the sequence
$$
x_0^1, x_0^2,\dots,x_0^s, x_1^1, \dots, x_1^s, x_2^1, \dots
$$
has the required properties.

Let $g=\k_Y(f)\in E(M)$. We construct by induction a sequence $\{z_n:n\in\o\}$ 
of points of $A$ such that:
\begin{enumerate}
\item if $k\le n$, then $d(z_n, a_k)=g(a_k)$;
\item $d(z_{n+1}, z_n)\le 2r_n$ for every $n\in \o$.
\end{enumerate}
Pick $z_0\in A$ so that $d(z_0, a_0)=g(a_0)$. This is possible since $A$ is
Urysohn. Suppose that the points $z_0, \dots, z_n$ have been constructed so
that the conditions 1 and 2 are satisfied. Consider two Kat\v etov
functions $f_n$ and $g_n$ on the set $A_{n+1}=\{a_k:k\le n+1\}$: let
$f_n(x)=d(z_n,x)$ for every $x\in A_{n+1}$, and let $g_n=g|_{A_{n+1}}$. By (1),
the functions $f_n$ and $g_n$ coincide on $A_n$, hence the distance between
them in the space $E(A_{n+1})$ is equal to
\begin{equation*}
\begin{split}
&|f_n(a_{n+1})-g_n(a_{n+1})|=
\sup\{|f_n(a_{n+1})-f_n(x)-g_n(a_{n+1})+g_n(x)|:x\in A_n\}\\
&\le\sup\{|f_n(a_{n+1})-f_n(x)|:x\in A_n\}+
\sup\{|g_n(a_{n+1})-g_n(x)|:x\in A_n\}\\
&\le 2d(a_{n+1}, A_n)=2r_n.
\end{split}
\end{equation*}

Let $X_n$ be the metric space $A_{n+1}\cup\{f_n\}$,
considered as a subspace of $E(A_{n+1})$. In virtue of Lemma~\ref{l:dist},
the map of $X_n$ to $A$ which leaves each point of $A_{n+1}$ fixed and sends
$f_n$ to $z_n$ is an isometric embedding. Since $A$ is Urysohn, this map can
be extended to an isometric embedding of $X_n\cup\{g_n\}$ to $A$. Let
$z_{n+1}$ be the image of $g_n$. Then $d(z_{n+1}, z_n)=d_{A_{n+1}}^E(g_n,f_n)
\le 2r_n$. In virtue of Lemma~\ref{l:dist}, for every $k\le n+1$ we have
$d(z_{n+1}, a_k)=g_n(a_k)=g(a_k)$. Thus the conditions~1 and~2 are
satisfied, and the construction is complete.

Since the series $\sum r_n$
converges, it follows from (2) that the sequence $\{z_n:n\in\o\}$ is
Cauchy and hence has a limit in the complete space $M$. Let
$z=\lim z_n$. By (1), we have $d(z,a_k)=g(a_k)$ for every $k\in \o$.
Since $Y$ is contained in the closure of the set $\{a_n:n\in \o\}$, it follows
that $d(z,y)=g(y)=f(y)$ for every $y\in Y$.
\end{proof}

\begin{proof}[Proof of Theorem~\ref{th:UrExt}]
Let $M$ be a metric space of diameter $\le1$, and let $M_\o$ be the
Urysohn extension of $M$ constructed above. Consider the completion $L$
of $M_\o$. Proposition~\ref{p:ThUr} implies that $L$ is Urysohn. 
Proposition~\ref{p:Katweight} shows that $w(L)=w(M)$. Finally, $M$ is
$g$-embedded in $M_\o$ (Proposition~\ref{p:propur}) and $M_\o$ is $g$-embedded
in $L$ (Proposition~\ref{p:compl}), so 
$M$ is $g$-embedded in $L$. Thus $L$ has the properties required by
Theorem~\ref{th:UrExt}.
\end{proof}

\section
{Graev metrics and $\o$-homogeneous extensions}
\label{s:graev}

In this section we prove the following:

\begin{thm}
\label{th:homogenext}
Every metric space can be $g$-embedded into
an $\o$-homogeneous metric space of the same weight and the same diameter.
\end{thm}

The proof is based on the construction of Graev metrics described in Section~2.
We apply this construction to metric spaces of relations. A {\it relation\/}
on a set $X$ is a subset of $X^2$. If $R$ and $S$ are relations on $X$, then
the composition $R\circ S$ (or simply $RS$) is defined by $R\circ S=
\{(x,y): \exists z((x,z)\in S \text{ and } (z,y)\in R)\}$. The inverse relation
$R\obr$ is defined by $R\obr=\{(x,y): (y,x)\in R\}$.
The set of all relations on a set $X$
is a semigroup with an involution:
the multiplication is given by the composition,
and the involution is given by the map $R\mapsto R\obr$. The unity
of this semigroup is the diagonal $\D$ of $X^2$.

We use the notation of Section~2. In particular, if $k$ is a non-expanding
function $\ge0$ 
on a metric space $(X,d)$, then $Gr(d,k)$ is the Graev pseudometric on
the free group $F(X)$. We consider the group $F(X)$ as a subset of the free
monoid $S(X)$ on the set $X\cup X\obr$.

\begin{proof}[Proof of Theorem~\ref{th:homogenext}]
Let $(M,d)$ be a metric space.
We first construct a $g$-embedding of
$M$ into a metric space $M^*$ such that $w(M^*)=w(M)$ and
every isometry between finite subsets of $M$ extends to an isometry of $M^*$.

For  every  isometry  $f:A\to  B$  between
finite  non-empty subsets of $M$ consider the graph $R=\{(a, f(a)):a\in A\}$ of
$f$, and let $\G$ be the set of all such graphs. Thus a non-empty finite subset
$R\sbs M^2$ is  an  element  of  $\G$  iff  for  any  two  pairs  $(x_1,y_1),\,
(x_2,y_2)\in  R$  we  have $d(x_1,x_2)=d(y_1,y_2)$. Equip $M^2$ with the metric
$d_2$ defined by  $d_2((x_1,y_1),(x_2,y_2))=  d(x_1,x_2)+d(y_1,y_2)$,  and  let
$d_H$  be  the  corresponding  Hausdorff metric on the set of finite subsets of
$M^2$. If $R$ and $S$ are two non-empty finite subsets of $M^2$ and  $a\ge  0$,
then  $d_H(R,S)\le  a$  iff  for every $p\in R$ there exists $q\in S$ such that
$d_2(p,q)\le a$, and for every  $p\in  S$  there  exists  $q\in  R$  such  that
$d_2(p,q)\le a$.

Let   $k$   be   the   non-expanding   function   on  $(\G,  d_H)$  defined  by
$k(R)=\max\{d(x,y): (x,y)\in R\}$. Let $G$ be the free group on $\G$,  equipped
with  the  Graev  pseudometric $D=Gr(d_H,k)$. To
avoid confusion of multiplication in $G$  with  composition  of  relations,  we
assign  to  each  $R\in  \G$  a  symbol  $t_R$, and consider elements of $G$ as
irreducible   words   of   the   form   $x_1^{\e_1}\dots   x_n^{\e_n}$,   where
$x_i=t_{R_i}$.  Similarly, we consider elements of the semigroup $S(\G)$ as
words of the same form.
Let  $\D=\{(x,x):x\in  M\}$  be  the  diagonal  of  $M^2$.
The set $\G'=\G\cup\{\emptyset\}\cup\{\D\}$ is a symmetrical subsemigroup
of the semigroup of all relations on $M$.
Let
$\Phi:G\to\G'$ be the map defined  by  the  following
rule: if $w=t_{R_1}^{\e_1}\dots t_{R_n}^{\e_n}\in G$ is a non-empty irreducible
word,  then  $\Phi(w)= R_1^{\e_1}\circ \dots \circ R_n^{\e_n}$. If $a, b\in M$,
then $(a,b)\in\Phi(w)$ iff there exists a chain $c_0=b, c_1, \dots,  c_n=a$  of
points  of  $M$ such that for every $i=1, \dots, n$ we have either $\e_i=1$ and
$(c_i, c_{i-1})\in R_i$ or $\e_i=-1$ and $(c_{i-1}, c_i)\in R_i$. For the empty
word $e_G\in G$ we put $\Phi(e_G)=\D$.

Note that the definition of $\Phi(w)$ makes sense also without the assumption
that the word $w$ be irreducible, so we can assume that $\Phi$ is defined on
the set $S(\G)$ of all words of the form
$t_{R_1}^{\e_1}\dots t_{R_n}^{\e_n}$. 
Recall that $w_1|w_2$ denotes the word obtained by
writing $w_2$ after $w_1$ (without cancelations).
We have $\Phi(w_1|w_2)=\Phi(w_1)\circ\Phi(w_2)$.

\begin{lemma}
\label{l:newphi}
If $w\in S(\G)$ and $u$ is the irreducible word equivalent to $w$,
then $\Phi(u)\supset \Phi(w)$. 
\end{lemma}

\begin{proof}
It suffices to prove that $\Phi(w')\supset \Phi(w)$
if $w'$ is obtained from $w$ by canceling one pair of
letters.  
Let $w=u|t_R^\e t_R^{-\e}|v$ and $w'=u|v$.
Since $R^\e$ is a functional relation, we have 
$R^\e\circ R^{-\e}\sbs\D$ and hence
$\Phi(w')=\Phi(u)\circ \Phi(v)=\Phi(u)\circ \D\circ\Phi(v)
\supset\Phi(u)\circ R^\e\circ R^{-\e}
\circ \Phi(v)=\Phi(w)$.
\end{proof}

For every $w\in G$ we have $\Phi(w\obr)=\Phi(w)\obr$. We claim that
$\Phi(w_1w_2)\supset \Phi(w_1)\circ \Phi(w_2)$ for every $w_1, w_2\in G$.
Indeed, the product $w_1w_2\in G$
is the irreducible word equivalent to $w_1|w_2$, therefore
$\Phi(w_1w_2)\supset
\Phi(w_1|w_2)=\Phi(w_1)\circ \Phi(w_2)$ by Lemma~\ref{l:newphi}.

For every $a, b\in M$ let $H_{a,b}\sbs G$ be the set of all $w\in G$ such  that
$(a,b)\in  \Phi(w)$.
We claim that $H_{a,b}\obr=H_{b,a}$ and $H_{b,c}H_{a,b}\sbs H_{a,c}$ for
every $a,b,c\in M$.
This follows from the  properties  of  $\Phi$  established in the preceding
paragraph. Indeed,
pick $w_1\in H_{b,c}$ and $w_2\in H_{a,b}$. Then $(a,b)\in \Phi(w_2)$ and
$(b,c)\in \Phi(w_1)$, hence $(a,c)\in \Phi(w_1)\circ \Phi(w_2)\sbs
\Phi(w_1w_2)$ and $w_1w_2\in H_{a,c}$. This proves
the inclusion $H_{b,c}H_{a,b}\sbs H_{a,c}$.
The equality $H_{a,b}\obr=H_{b,a}$
is proved similarly.

Note that $t_R\in H_{a,b}$ if and only if $(a,b)\in R$, since $\Phi(t_R)=R$.
Note also that $e_G\in H_{a,b}$ if and only if $a=b$, since $\Phi(e_G)=\D$.

Consider the following equivalence relation  $\sim$  on  $G\times  M$:  a  pair
$(g,a)$  is  equivalent  to  a  pair  $(h,b)$  iff  $h\obr g\in H_{a,b}$. Since
$e_G\in H_{a,a}$,
$H_{a,b}\obr=H_{b,a}$  and  $H_{b,c}H_{a,b}\sbs H_{a,c}$ for all $a,b,c\in M$,
the relation $\sim$ is reflexive, symmetric and transitive and thus is indeed 
an equivalence relation.
Let $L$ be the quotient set $G\times M/\sim$. The group
$G$  acts on $G\times M$ by the rule $g\cdot (h,a)=(gh,a)$. The relation $\sim$
is invariant under this action, so there is a uniquely defined left  action  of
$G$  on  $L$ which makes the canonical map $G\times M\to L$ into a morphism
of $G$-sets.
Let $i:M\to L$ be the map which sends each point  $a\in
M$ to the class of the pair $(1,a)$. If $a\ne b$, then the pairs $(1,a)$ and
$(1,b)$ are not equivalent, since $e_G\notin H_{a,b}$. The map $i$ is therefore
injective, and we can consider $M$ as a subspace of $L$, identifying $M$ with
$i(M)$. Every $x\in L$ can be written in the form $x=g\cdot a$
(or simply $x=ga$), where $g\in G$
and $a\in M$. 

Let $a, b\in M$. The set of all $g\in G$ such that $ga=b$ is equal to
$H_{a,b}$. If $R\in \G$ is a relation containing the pair  $(a,b)$,
then  $t_R\in  H_{a,b}$  and hence $t_Ra=b$. It follows that the action of $G$
on $L$ is transitive. Moreover, for every isometry $f:A\to B$ between finite
subsets of $M$ there exists $g\in G$ such that the self-map $x\to gx$ of $L$
extends $f$. Indeed, if $R\in \G$ is the graph of $f$, then 
$t_R\in H_{a,f(a)}$ and hence $t_Ra=f(a)$ for every $a\in A$.
Thus $g=t_R$ has the required property.

We now define a $G$-invariant pseudometric $\nu$
on $L$ which extends the metric $d$ on $M$. 
Let $p$ be the Graev seminorm on $G$ corresponding to the pseudometric
$D=Gr(d_H,k)$. We have $p(w)=D(w,e_G)$ for every $w\in G$.
For every $x,y\in L$ let
$$
\nu(x,y)=\inf \{p(g): g\in G,\  gx=y\}.
$$
Then $\nu$ is a pseudometric on $L$. Since the seminorm $p$ is invariant
under inner automorphisms, the pseudometric $\nu$ is $G$-invariant. Indeed,
for $x,y\in L$ and $h\in G$ we have $\nu(hx,hy)=\inf\{p(g):ghx=hy\}=
\inf\{p(h\obr gh):h\obr ghx=y\}=\inf\{p(g'):g'x=y\}=\nu(x,y)$.
We claim that $\nu$ extends the metric $d$ on $M$:
$d(a,b)=\nu(a,b)$  for  every $a, b\in M$.
Since for $w\in G$ the condition  $wa=b$  is  equivalent  to $w\in H_{a,b}$, 
we have $\nu(a, b) = \inf\{p(w):w\in  H_{a,b}\}$. If $R=\{(a,b)\}$,  
then $t_R\in H_{a,b}$   and
$p(t_R)=k(R)=d(a,b)$.  It follows that $\nu(a,b) \le d(a,b)$. 
It remains to prove
the opposite inequality, which is equivalent to the following assertion:

\begin{lemma}
\label{l:udine}
If $a, b\in M$ and $w\in H_{a,b}$, then $p(w)\ge d(a,b)$.
\end{lemma}

\begin{proof} 
Let $w= t_{R_1}^{\e_1}\dots t_{R_n}^{\e_n}$. We argue by induction
on $n$, the length of $w$.  
If $n=0$, then $w=e_G$, and we noted that $e_G\in H_{a,b}$ implies $a=b$.
If  $n=1$,  then  $w=t_R^\e$ and $p(w)=k(R)$.
Since $w\in H_{a,b}$, the relation  $R$  contains either
$(a,b)$  or  $(b,a)$  and hence $p(w)=k(R)\ge d(a,b)$. Assume that $n>1$. It
suffices to show that there exists $u\in H_{a,b}$ of  length  $<n$  such  that
$p(u)\le p(w)$.

We use the construction of the Graev seminorm $p$ described in Section~2.
Let  $E$  be a $w$-pairing for which $p(w)$ is attained. In other
words, $E$ is a disjoint  system  of  two-element  subsets  of  the set
$J=\{1,\dots,n\}$  such  that  for  the  Graev sum
$$ s_E =\sum \{d_H(R_i,R_j):
\{i,j\}\in E,\,i<j\} +
 \sum \{k(R_i):i\in J\setminus \cup E\}
$$
we have $p(w)=s_E$. Considering the pair $(i,j)\in E$
with the least possible value of $|i-j|$
(``the shortest arc"), we see  that  at  least  one  of  the
following  three  cases must occur:
\begin{enumerate}
\item there exists an $i$ such that
$\{i,i+1\}\in E$;
\item there exists an $i$ such that  $\{i,i+2\}\in  E$  and
$i+1\in  J\setminus  \cup  E$;
\item  there  exists an $i$ such that $i,i+1\in
J\setminus \cup E$.
\end{enumerate}
In cases (1)  or  (3)  we  replace  the  subword
$t_{R_i}^{\e_i}t_{R_{i+1}}^{\e_{i+1}}$  of  $w$  by  the  letter  $t_S$,  where
$S=R_i^{\e_i}\circ R_{i+1}^{\e_{i+1}}$. In case  (2)  we  replace  the  subword
$t_{R_i}^{\e_i}t_{R_{i+1}}^{\e_{i+1}}t_{R_{i+2}}^{\e_{i+2}}$   of  $w$  by  the
letter  $t_S$,  where     
$S=R_i^{\e_i}\circ R_{i+1}^{\e_{i+1}}\circ R_{i+2}^{\e_{i+2}}$.
In  all  cases we get a word $w'$ of length $<n$. To justify the usage of the
symbol $t_S$, we must show that $S\in \G$, which reduces to the fact that 
$S\ne\emptyset$. Had $S$ been empty, the same would have been true for
$\Phi(w)=R_1^{\e_1}\circ\dots\circ R_n^{\e_n}$. On the other hand, since
$w\in H_{a,b}$, we have $(a,b)\in \Phi(w)\ne\emptyset$.

Let $u\in G$ be
the irreducible word equivalent to $w'$. The length of $u$ 
is less than $n$.
We show
that $u\in H_{a,b}$ and $p(u)\le p(w)$.

By Lemma~\ref{l:newphi} we have
$\Phi(w')\sbs\Phi(u)$. Plainly
$\Phi(w)=R_1^{\e_1}\circ\dots\circ R_n^{\e_n}=\Phi(w')$.
Since $w\in H_{a,b}$, we have  $(a,b)\in  \Phi(w)$.  Thus
$(a,b)\in\Phi(w)=\Phi(w')\sbs\Phi (u)$ and $u\in H_{a,b}$, as required.

We  prove that $p(u)\le p(w)$. As in Section~2, we define $p(w')$ even if
the word $w'$ is reducible, and we have $p(u)=p(w')=\inf s_{E'}$, where
$E'$ runs over the set of all $w'$-pairings.
The $w$-pairing $E$ in an obvious way yields 
a $w'$-pairing $E'$, which coincides with $E$ outside the changed part
of $w$ and leaves the new letter $t_S$  unpaired.  The  Graev  sums  $s_E$  and
$s_{E'}$  differ  only  by  the  term  $k(S)$ in the sum $s_{E'}$ and the terms
$d_H(R_i,R_{i+1})$ (case 1)  or  $d_H(R_i,R_{i+2})  +k(R_{i+1})$  (case  2)  or
$k(R_i)+k(R_{i+1})$ (case 3) in the sum $s_E$. According to 
Lemma~\ref{l:udikrana} below, we
have  $s_{E'}\le  s_E$.  
Thus $p(u)=p(w')\le s_{E'}\le s_E = p(w)$. 
\end{proof}

\begin{lemma}
\label{l:udikrana}
Let $\e,\delta\in\{-1,1\}$.
\begin{enumerate}
\item If $R_1,R_2\in \G$ and $S=R_1^\e R_2^{-\e}$ is non-empty, then
$k(S)\le d_H(R_1,R_2)$;
\item if $R_1,R_2,R_3\in\G$ and $S=R_1^\e R_2^\delta R_3^{-\e}$ is non-empty,
then $k(S)\le d_H(R_1,R_3)+k(R_2)$;
\item if $R_1,R_2\in\G$ and $S=R_1^\e R_2^\delta$ is non-empty, then
$k(S)\le k(R_1)+k(R_2)$.
\end{enumerate}
\end{lemma}

\begin{proof}
Since $k(R)=k(R\obr)$ and $d_H (R,T)=d_H(R\obr,T\obr)$ for every $R,T\in\G$, we
may  assume  that  $\e=\delta=1$.  
Pick $(a,b)\in S$ so that $k(S)=d(a,b)$.
Case (1) follows from (2) 
(take for $R_2$ in (2) a sufficiently large finite part of $\D$), 
so let us consider case (2).
There exist $x,y\in M$  such
that  $(a,x)\in  R_3\obr$,
$(x,y)\in R_2$ and $(y,b)\in R_1$. Since $(x,a)\in R_3$,
there exists a pair $(u,v)\in R_1$ such that $d(a,v)+d(u,x)\le d_H(R_1,  R_3)$.
The  relation  $R_1$,  being  an  element  of  $\G$,  is the graph of a partial
isometry,  so  from  $(y,b)\in  R_1$  and  $(u,v)\in  R_1$  it   follows   that
$d(v,b)=d(u,y)$. Note that $d(x,y)\le k(R_2)$.
Thus we have $k(S)=d(a,b)\le d(a,v)+d(v,b)=d(a,v)+d(u,y)\le
d(a,v)+d(u,x)+d(x,y)\le d_H(R_1,R_3)+k(R_2)$, as required.
Case  (3)  is  easy: there exists a point $c\in M$ such that $(a,c)\in R_2$ and
$(c,b)\in  R_1$,   hence   $k(S)=d(a,b)\le   d(a,c)+d(c,b)\le   k(R_2)+k(R_1)$.
\end{proof}

We have thus proved that the pseudometric $\nu$ on $L$ extends the metric
$d$ on $M$.
Let $(M^*,d^*)$ be the metric space
associated with the pseudometric space $(L,\nu)$.
The metric space 
$(M,d)$ can be naturally identified with a subspace of $(M^*,d^*)$. We
show that $M$ is $g$-embedded in $M^*$.

In virtue of the functorial nature of the construction of $M^*$, every
isometry $\f$ of $M$ naturally extends to an isometry $\f^*$ of $M^*$.
The map $\f\mapsto\f^*$ from $\Is(M)$ to $\Is(M^*)$ is a homomorphism 
of groups. We claim that this homomorphism is continuous. This follows
from the fact that at each step of our construction new spaces are obtained
from the old ones via functors ``with finite support": every element of 
$\G$ is a {\it finite\/} relation on $M$, and every word $w\in G$ involves
only finitely many elements of $\G$.
Given an isometry $\f\in \Is(M)$, 
the isometry $\f^*\in \Is(M^*)$ can be obtained step by step
in the following way. 
First we consider the isometry $\f_1$ of the metric space $(\G, d_H)$
corresponding to $\f$; the isometry $\f_1$ preserves the function $k$ on $\G$
and gives rise to the automorphism
$\f_2$ of the group $G=F(\G)$ which preserves the Graev pseudometric $D$;
then we get the isometry $\f_3$ of $L$ which maps the class of each pair
$(g,x)$ ($g\in G$, $x\in M$) to the class of the pair $(\f_2(g), \f(x))$;
finally we get the isometry $\f_4=\f^*$ of $M^*$. We show step by step that 
$\f_i$ depends continuously on $\f$. For $i=1$ this is straightforward:
use the fact that $\G$ consists of {\it finite\/} subsets of $M^2$. For $i=2$
apply Lemma~\ref{l:Gr-metr-isom} with $X=\G$. Let us consider the case
$i=3$. Pick a point $x=ga\in L$ ($g\in G$, $a\in M$). It suffices to check that
$\nu(\f_3(x),x)$ is small if $\f$ is close to the identity. We have
$\nu(\f_3(x),x)=\nu(\f_2(g)\f(a), ga)\le \nu(\f_2(g)\f(a), g\f(a))+
\nu(g\f(a), ga)=\nu(g\obr\f_2(g)\f(a),\f(a))+\nu(\f(a),a)$.
By the definition of $\nu$,
the first term of the last sum does not exceed $p(g\obr\f_2(g))=D(\f_2(g),g)$
and hence is arbitrarily small if $\f$ is close enough to the identity.
The same is true for second term, and we are done. Finally, $\f_4$ is the
image of $\f_3$ under the natural morphism $\Is(L)\to\Is(M^*)$, and the case
$i=4$ follows.

We have thus proved that $M$ is $g$-embedded in $M^*$. We saw that each 
isometry between finite subsets of $M$ extends to an isometry of $L$ and hence
also to an isometry of $M^*$. It is easy to see that $w(M^*)=w(M)$.
If the diameter $C$ of $M$ is finite, replace the metric $d^*$ of $M^*$
by 
$\inf(d^*,C)$. This operation can make the group $\Is(M^*)$ only larger, 
and the diameter of $M^*$ becomes equal to that of $M$.

To finish the proof of Theorem~\ref{th:homogenext}, iterate the construction
of $M^*$. We get an increasing chain 
$M_0=M\sbs M_1=M^*\sbs M_2=M_1^* \sbs\dots$ of
metric spaces such that each $M_n$ is $g$-embedded in $M_{n+1}$, every isometry
between finite subsets of $M_n$ extends to an isometry of $M_{n+1}$,
$w(M_n)=w(M)$ and $\diam M_n=\diam M$, $n=0,1,\dots$.
Consider the space $M_\o=\bigcup_{n\in\o}M_n$.
We have $w(M_\o)=w(M)$ and $\diam M_\o = \diam M$.
In virtue of Proposition~\ref{p:gchain}, each $M_n$ is $g$-embedded in $M_\o$.
Since every finite subset of $M_\o$ is contained in some $M_n$, it is clear
that $M_\o$ is $\o$-homogeneous.
\end{proof}

{\it Remarks.} 

1. If $a,b\in M$ are distinct 
and $S=\{(b,b)\}$, the pairs $(1,a)$ and $(t_S,a)$ represent distinct points
of $L$ that have the same image in $M^*$.
Early versions of this paper contained the false statement
that $\nu$ itself is a metric and $M^*=L$.  I am indebted to the referee for
catching this error. 

2. The referee raised the question whether the methods of this section could
be used to prove the following result by S. Solecki \cite{Sol}
and A.M. Vershik \cite{Versh07} that extends an earlier result by Hrushovski:
for every finite metric space $A$ there exists another finite metric space
$A^*$ containing $A$ such that all partial isometries%
\footnote
{A {\em partial isometry} of $A$ is an isometry between two subsets of $A$.}
 of $A$ extend to isometries
of $A^*$. I do not know the answer. A partial answer is provided by
Pestov's paper \cite{Pes07} where the Hrushovski--Solecki--Vershik theorem
is proved with the aid of pseudometrics on groups, 
and the notion of a residually
finite group is used to construct isometric embeddings of finite metric
spaces into finite metric groups. A similar technique
was used in \cite{PU}.

\section
{Proof of Theorem~\ref{th:1.4}}
\label{s:proof1.4}

In this section we prove Theorem~\ref{th:1.4}.

\begin{thm}
\label{th:pred1.4}
Let $M$ be a metric space of diameter $\le1$. There exists a complete
$\o$-homogeneous Urysohn metric space $L$ containing $M$ as a subspace 
such that $w(L)=w(M)$ and $M$ is $g$-embedded in $L$.
\end{thm}

\begin{proof}
Consider two cases.

Case 1: $M$ is separable. According to Theorem~\ref{th:UrExt}, there exists a
complete separable Urysohn space $L$ such that $M$ is a $g$-embedded 
subspace of $L$. According to Proposition~\ref{p:newequiv}, 
$L=\U_1$ is $\o$-homogeneous.
%
%

Case 2: $M$ is not separable. Let $\tau=w(M)$.
Applying in turn Theorem~\ref{th:UrExt} and
Theorem~\ref{th:homogenext}, construct an increasing continuous chain
$\{M_\a:\a\le\o_1\}$ of metric spaces of weight $\tau$ and diameter $\le1$
such that $M_0=M$, each $M_\a$ is $g$-embedded in $M_{\a+1}$ ($\a<\o_1$),
and $M_{\a+1}$ is complete Urysohn for $\a$ even and $\o$-homogeneous for 
$\a$ odd. Let $L=M_{\o_1}=\bigcup_{\a<\o_1}M_\a$.
Proposition~\ref{p:gchain} implies that each $M_\a$ is $g$-embedded
in $L$. The space $L$ is Urysohn, being the union of the increasing
chain $\{M_{2\a+1}:\a<\o_1\}$ of Urysohn spaces. For similar reasons the
space $L$ is $\o$-homogeneous. Finally, since
every countable subset of $L$ is contained in some $M_\a$, $\a<\o_1$,
and all spaces $M_{2\a+1}$ are complete,
every Cauchy sequence in $L$ converges, which means that $L$ is complete.
Thus $L$ has the properties required by Theorem~\ref{th:pred1.4}.
\end{proof}

\begin{proof}[Proof of Theorem~\ref{th:1.4}]
Let $G$ be a topological group. According to Theorem~\ref{th:GrIsom}, there
exists a metric space $(M,d)$ such that $w(M)=w(G)$ and $G$ is isomorphic to
a subgroup of $\Is(M)$. We may assume that $M$ has diameter $\le1$: otherwise
replace the metric $d$ by 
$\inf(d,1)$. Theorem~\ref{th:pred1.4} implies that
there exists a complete $\o$-homogeneous Urysohn metric space $L$ such that 
$w(L)=w(M)$ and $\Is(M)$ is isomorphic to a subgroup of $\Is(L)$. Then
$w(L)=w(G)$ and $G$ is isomorphic to a subgroup of $\Is(L)$, as required.
\end{proof}

\section{Semigroups of bi-Kat\v etov functions}
\label{s:semi}

Let $(M,d)$ be a complete metric space of diameter $\le1$.

\begin{defin}
\label{def:bikat}
A function $f:M\ti M\to I=[0,1]$ is {\it bi-Kat\v etov\/} if for every $x\in M$
the functions $f(x, \cdot)$ and $f(\cdot, x)$ on $M$ are Kat\v etov 
(see Section~\ref{s:katconstr}).
\end{defin}

Thus a function $f:M^2\to I$ is bi-Kat\v etov if and only if for every
$x,y,z\in M$ we have
\begin{gather*}
|f(x,y)-f(x,z)|\le d(y,z)\le f(x,y)+f(x,z),\\
|f(y,x)-f(z,x)|\le d(y,z)\le f(y,x)+f(z,x).
\end{gather*}
Let $\Th$ be the compact space of all bi-Kat\v etov
functions on $M^2$, equipped
with the topology of pointwise convergence.
In the next section we shall prove that the Roelcke completion
of the group $\Is(M)$ can be identified with $\Th$, provided that 
the complete metric
space $M$ is Urysohn and $\o$-homogeneous. In the present section
we study the structure of an ordered semigroup with an involution on $\Th$.

Recall that we defined in Section~1 an associative operation $\bull$ on the set
$S=I^{M\times M}$. If $f,g\in S$ and $x,y\in M$, then
$$
f\bull g(x,y)=\inf\{f(x,z)\uplus g(z,y):z\in M\}.
$$
The involution $f\mapsto f^*$ on $S$ is defined by $f^*(x,y)=f(y,x)$. 
Every idempotent in $S$ satisfies the triangle inequality.
If $f\in S$ is zero on the diagonal of $M^2$, then $f$ is an idempotent in $S$
if and only if $f$ satisfies the triangle inequality.
A function 
$f\in S$ is a pseudometric on $X$ if and only if $f$ is zero on the diagonal
and $f$ is a symmetrical idempotent. In particular, we have $d=d^*=d\bull d$.

The semigroup $S$ has a natural partial order: for $p,q\in S$ the inequality
$p\le q$ means that $p(x,y)\le q(x,y)$ for all $x,y\in M$. This partial
order is compatible with the semigroup structure: if $p_1\le p_2$ and 
$q_1\le q_2$, then $p_1\bull q_1\le p_2\bull q_2$.

It is clear that the set $\Th$ of all bi-Kat\v etov functions is closed under
the involution. It is easy to verify
that $\Th$ also is closed under the operation $\bull$. This fact also can be 
deduced from the following proposition:

\begin{prop}
\label{p:charbikat}
A function $f:M^2\to I$ is bi-Kat\v etov if and only if
$$
f\bull d=d\bull f=f, \quad f^*\bull f\ge d,\quad f\bull f^*\ge d,
$$
where $d$ is the metric on $M$.
\end{prop}

\begin{proof}
The condition $f\bull d=f$ (respectively, $d\bull f=f$) holds if and only if
the function $f(x,\cdot)$ (respectively, $f(\cdot,x)$) is non-expanding for
every $x\in X$. The condition $f^*\bull f\ge d$
(respectively, $f\bull f^*\ge d$) holds if and only if
$d(y,z)\le f(x,y)+f(x,z)$ (respectively, $d(y,z)\le f(y,x)+f(z,x)$) for all
$x,y,z\in X$.
\end{proof}

Let $S$ be any ordered semigroup with an involution, and let $d\in S$ be
a symmetrical idempotent. The set $S_d$ of all $x\in S$ such that
$$
xd=dx=x, \quad x^*x\ge d, \quad xx^*\ge d
$$
is closed under the multiplication and under the involution and can 
be considered as a semigroup with the unity $d$. Indeed, we have $d\in S_d$
since $d=d^*=d^2$, and it is clear that $d$ is the unity of $S_d$. If 
$x,y\in S_d$, then $xyd=xy=dxy$ and $(xy)^*xy=y^*x^*xy\ge y^*dy=y^*y\ge d$;
similarly, $xy(xy)^*\ge d$ and hence $xy\in S_d$. Thus $S_d$ is a semigroup.
If $x\in S_d$, then $x^*d=x^*d^*=(dx)^*=x^*$ and similarly $dx^*=x^*$. It
follows that $S_d$ is symmetrical.

The arguments of the preceding paragraph and Proposition~\ref{p:charbikat}
show that $\Th$ is a semigroup with the unity $d$.
In general,  the operation $(f,g)\mapsto f\bull g$
need not be continuous (not even continuous on the left or on the right).

\begin{prop}
\label{p:maximidem}
Let $S$ be a closed subsemigroup of $\Th$, and let $T$ be the set of all 
$f\in S$ such that $f\ge d$. If $T\ne\emptyset$,
then $T$ has a greatest element $p$, and
$p$ is an idempotent.
\end{prop}

\begin{proof}
We claim that every non-empty closed subset of $\Th$ has a maximal element.
Indeed, if $C$ is a non-empty
linearly ordered subset of $\Th$, then $C$ has a least upper bound $b$ in $\Th$,
and $b$ belongs to the closure of $C$. Thus our claim follows from
Zorn's lemma. 

The set $T$ is a closed subsemigroup of $\Th$.
Let $p$ be a maximal element of $T$.
For every $q\in T$ we have $p\bull q\ge p\bull d=p$, whence $p\bull q=p$.
It follows that $p$ is idempotent and that 
$p=p\bull q\ge d\bull q=q$. Thus $p$ is the greatest element of $T$.
\end{proof}

We now describe all idempotents in $\Th$ which are $\ge d$. For every closed
non-empty
subset $F$ of $M$ let $b_F\in \Th$ be the bi-Kat\v etov function defined by
$b_F(x,y)=\inf\{d(x,z)\uplus d(z,y):z\in F\}$. If $F=\emptyset$, let $b_F=1$,
that is the function on $M^2$ which is identically equal to 1. (Note that~1
is {\it not\/} the unity of $\Th$;
on the contrary, $f\bull 1=1\bull f=1$ for every $f\in I^{M\ti M}$,
so 1 might be called a zero element of $\Theta$.)

\begin{prop}
\label{p:idemp}
If $F$ is a closed subset of $M$, then $b_F$ is an idempotent $\ge d$
in $\Th$, and 
every idempotent $\ge d$ in $\Th$ is equal to $b_F$ for some closed $F\sbs M$.
\end{prop}

\begin{proof}
Let $F$ be a closed subset of $M$. It is clear that $b_F\ge d$. 
If $F\ne\emptyset$, then $b_F\bull b_F(x,y)=
\inf\{d(x,z_1)\uplus d(z_1,u)\uplus d(u,z_2)\uplus d(z_2,y):
u\in M,\,z_1,z_2\in F\}=\inf\{d(x,z)\uplus d(z,y):z\in F\}
=b_F(x,y)$
for every $x,y\in M$. Thus $b_F$ is an idempotent. The same is obviously true
if $F=\emptyset$.

Conversely, let $p$ be an idempotent in $\Th$ such that $p\ge d$. 
Let $F=\{x\in M:p(x,x)=0\}$. The function $p:M^2\to I$, being non-expanding in
each argument, is continuous, hence
$F$ is closed in $M$. We claim that $p=b_F$.

We first show that $p\le b_F$. This is evident if $F=\emptyset$, so assume
that $F\ne\emptyset$. For every $x,y,z\in M$ we have $p(x,y)\le d(x,z)+p(z,y)
\le d(x,z)+d(z,y)+p(z,z)$, since the functions $p(\cdot,y)$ and $p(z,\cdot)$
are non-expanding. It follows that $p(x,y)\le 
\inf(\{d(x,z)+d(z,y)+p(z,z):z\in F\}\cup\{1\})=b_F(x,y)$.

We prove that $b_F\le p$. Fix $x,y\in M$. We must show that
$b_F(x,y)\le p(x,y)$. This is evident if $p(x,y)=1$, so assume that $p(x,y)<1$.
Fix $\e>0$ so that $p(x,y)+\e<1$. Since $p\bull p=p$, for every $u,v\in M$ we 
have $p(u,v)=\inf(\{p(u,z)+p(z,v):z\in M\}\cup \{1\})$. Hence we can construct
by induction a sequence of points $z_1,z_2,\dots$ in $M$ such that 
\begin{align*}
p(x,z_1)+p(z_1,y)&<p(x,y)+\e/2;\\
p(z_1,z_2)+p(z_2,y)&<p(z_1,y)+\e/4;\\
p(z_2,z_3)+p(z_3,y)&<p(z_2,y)+\e/8;\\
..................&.................\\
p(z_n,z_{n+1})+p(z_{n+1}, y)&<p(z_n,y)+\e/2^{n+1}\\
..................&.................
\end{align*}
Adding the first $n$ inequalities, we get 
$$
p(x,z_1)+\sum_{i=1}^{n-1}p(z_i,z_{i+1})+p(z_n,y)<p(x,y)+\e. \leqno{\rm (I)}
$$ 
It follows
that the series $\sum_{i=1}^\infty p(z_i,z_{i+1})$ converges. Since $d\le p$,
the series $\sum_{i=1}^\infty d(z_i,z_{i+1})$ also converges. This implies
that the sequence $z_1,z_2,\dots$ is Cauchy and hence has a limit in $M$.
Let $z=\lim z_i$. Since the series $\sum_{i=1}^\infty p(z_i,z_{i+1})$ 
converges, we have $\lim p(z_i,z_{i+1})=0$ and hence $p(z,z)=0$. Thus $z\in F$.

Since 
$p$ is an idempotent, it 
satisfies the triangle inequality: $p(u,v)\le p(u,w)+p(w,v)$ for all $u,v,w\in M$.
The inequality (I) therefore implies that $p(x,z_n)+p(z_n,y)<p(x,y)+\e$ for
every $n$. Passing to the limit, we get
$p(x,z)+p(z,y)\le p(x,y)+\e$. 
Thus $b_F(x,y) \le d(x,z)+d(z,y)\le p(x,z)+p(z,y)\le p(x,y)+\e$.
Since $\e$ was arbitrary, 
it follows that $b_F(x,y)\le p(x,y)$.
\end{proof}

{\it Remark}. We shall see later in this section that the
elements of $\Th$ (= bi-Kat\v etov functions on $M^2$) admit a geometric
interpretation: they correspond to metric spaces covered by two isometric
copies of $M$. If $F$ is a closed subset of $M$, the function $b_F$ considered
above corresponds to the amalgam of two copies of $M$ with the copies of $F$
amalgamated. This description, together with the geometric description of the operation
$\bullet$ on $\Th$ provided in the last paragraph of this section, makes it obvious
that each $b_F$ is an idempotent. 

\smallskip

Let $G=\Is(M)$.
For every isometry $\f\in G$ let $i(\f)\in \Th$ be the bi-Kat\v etov function
defined by $i(\f)(x,y)=d(x,\f(y))$. It is easy to check that the map 
$i:G\to \Th$ is a homeomorphic embedding.
We claim that the embedding $i:G\to \Th$ is a morphism of monoids
with
an involution. This means that $i(e_G)=d$, $i(\f\obr)=i(\f)^*$ and
$i(\f\psi)=i(\f)\bull i(\psi)$ for all $\f,\psi\in G$. The first equality
is obvious. For the second, note that 
$i(\f\obr)(x,y)=d(x,\f\obr(y))=d(y,\f(x))=i(\f)(y,x)=i(\f)^*(x,y)$. 
For the third, note that $i(\f\psi)(x,y)=d(x,\f\psi(y))=
\inf\{d(x,\f(z))+d(\f(z),\f\psi(y)): z\in M\}=
\inf\{d(x,\f(z))+d(z,\psi(y)): z\in M\}=
\inf\{i(\f)(x,z)+i(\psi)(z,y):z\in M\}=i(\f)\bull i(\psi)(x,y)$.

Thus we can identify $G$ with a subgroup of $\Th$. There are natural
left and right actions of $G$ on $\Theta$, defined by $(g,p)\mapsto
g\bull p$ and $(g,p)\mapsto p\bull g$ ($g\in G,\ p\in \Theta$), respectively.

\begin{prop}
\label{p:contshift}
The maps $(g,p)\mapsto g\bull p$ and $(g,p)\mapsto p\bull g$
from $G\ti \Theta$ to $\Theta$ 
are continuous. If $p\in\Th$ and
$x,y\in M$, then $g\bull p(x,y)=p(g\obr(x),y)$ and
$p\bull g(x,y)=p(x,g(y))$.
\end{prop}

\begin{proof}
We have
$g\bull p(x,y)=\inf\{d(x,g(z))\uplus p(z,y):z\in M\}$. Taking $z=g\obr(x)$,
we see that the right side is $\le p(g\obr(x),y)$. On the other hand, for
every $z\in M$ we have $d(x,g(z))+p(z,y)=d(g\obr(x),z)+p(z,y)\ge 
p(g\obr(x),y)$, whence the opposite inequality. The continuity of 
the left action easily follows
from the explicit formula that we have just proved. The argument
for the right action is similar. 
\end{proof}

Let us show that all invertible elements of $\Th$ are in $i(G)$.

It will be useful to
establish a one-to-one correspondence between elements of $\Th$ and other
objects which we call $M$-triples. Let $s=(h_1,h_2,L)$ be a triple such
that $L$ is a metric space of diameter $\le1$, $h_i: M\to L$ is an isometric
embedding ($i=1,2$) and $L=h_1(M)\cup h_2(M)$. We say that $s$ is an
{\it $M$-triple\/}. Two $M$-triples $(h_1,h_2,L)$ and
$(h_1',h_2',L')$ are {\it isomorphic\/} if there exists an isometry
$g:L\to L'$ such that $h_i'=gh_i$, $i=1,2$.

Given an $M$-triple $s=(h_1,h_2,L)$, let $f_s\in \Th$ be the
bi-Kat\v etov function
defined by $f_s(x,y)=\rho_L(h_1(x),h_2(y))$, where $\rho_L$ is the metric
on $L$. It is easy to verify that we get in this way
a one-to-one correspondence between $\Th$ and the set
of classes of isomorphic $M$-triples. The subset $i(G)$ of $\Th$
corresponds to the set of classes of triples $s=(h_1,h_2,L)$ such that
$h_1(M)=h_2(M)=L$. Indeed, if $\f\in G$, then for the $M$-triple
$s=(\text{id}_M, \f, M)$ we have $f_s=i(\f)$. Conversely, every
$M$-triple $s=(h_1,h_2,L)$ such that 
$h_1(M)=h_2(M)=L$ is isomorphic to the triple
$(\text{id}_M, \f, M)$, where $\f=h_1\obr h_2$ is an isometry of $M$.
Thus $s$ corresponds to $\f\in G$.

\begin{prop}
\label{p:invert}
The set of invertible elements of $\Th$ coincides with $i(G)$.
\end{prop}

\begin{proof}
Let $f\in \Th$ be invertible. Let $s=(h_1,h_2,L)$ be an $M$-triple
corresponding to $f$. This means that $(L,\rho)$ is a metric
space, $h_1$ and $h_2$ are distance-preserving maps from $M$ to $L$, 
$L=h_1(M)\cup h_2(M)$ and $f(x,y)=\rho(h_1(x),h_2(y))$ for all $x,y\in M$.
We saw that elements of $G$ correspond to triples $s$ satisfying the condition
$h_1(M)=h_2(M)=L$. Thus we must verify this condition.

Let $g$ be the inverse of $f$. Then $f\bull g=g\bull f=d$. For every $x\in M$
we have $\inf\{f(x,y)+g(y,x):y\in M\}=f\bull g(x,x)=d(x,x)=0$ and hence
$\rho(h_1(x),h_2(M))=\inf\{f(x,y):y\in M\}=0$. This means that $h_1(x)$
belongs to the closure of $h_2(M)$ in $L$.
Since $M$ is complete and 
$h_2$ is an isometric embedding, $h_2(M)$ is closed in $L$. It follows that
$h_1(x)\in h_2(M)$. Since $x\in M$ was arbitrary, we have $h_1(M)\sbs h_2(M)$.
Similarly, $h_2(M)\sbs h_1(M)$ and therefore $h_1(M)=h_2(M)=L$.
\end{proof}

The operation $\bull$ has the following description in terms of 
$M$-triples. Let $p,q\in\Th$. There exists a quadruple $s=(h_1,h_2,h_3,L)$
such that $(L,\rho)$
is a metric space of diameter $\le 1$, $L=L_1\cup L_2\cup L_3$,
$h_i:M\to L_i$ is an isometry ($i=1,2,3$), 
$(h_1,h_2, L_1\cup L_2)$ is an $M$-triple corresponding to $p$ and
$(h_2,h_3, L_2\cup L_3)$ is an $M$-triple corresponding to $q$.
The bi-Kat\v etov function $f$ corresponding to the $M$-triple
$(h_1,h_3,L_1\cup L_3)$ depends on $s$, and the largest function $f$
over all quadruples $s$ such as above is equal to $p\bull q$. 
Indeed, we have $f(x,y)=\rho(h_1(x),h_3(y))\le 
\inf\{\rho(h_1(x),h_2(z))\uplus \rho(h_2(z), h_3(y)):z\in M\}=
\inf\{p(x,z)\uplus q(z,y):z\in M\}=p\bull q(x,y)$. To see that the function
$p\bull q$ can be attained, consider two disjoint copies $M'$ and $M''$ of $M$.
For $x\in M$ denote by $x'$ the copy of $x$ in $M'$, and use similar notation
for $M''$.
Let $\rho$ be the pseudometric on $X=M\cup M'\cup M''$ defined by 
$\rho(x,y)=\rho(x',y')=\rho(x'',y'')=d(x,y)$,
$\rho(x,y')=p(x,y)$, $\rho(x',y'')=q(x,y)$ and $\rho(x,y'')=p\bull q(x,y)$.
The triangle inequality for $\rho$ is easily verified. 
(The space $X$ is the amalgam 
(in the class of spaces of diameter $\le 1$)
of the spaces $M\cup M'$ and $M'\cup M''$
with the subspace $M'$ amalgamated, see \cite{UryHilb}
for a definition.)
Let $L$ be the metric
space associated with the pseudometric space $(X,\rho)$. Let $L_1,L_2,L_3$ be 
the images of $M$, $M'$, $M''$ in $L$, respectively. Let $h_i:M\to L_i$ be
the obvious isometry, $i=1,2,3$. The quadruple $s=(h_1,h_2,h_3,L)$ has the
properties considered above, and the bi-Kat\v etov function corresponding to
the $M$-triple $(h_1,h_3,L_1\cup L_3)$ is equal to $p\bull q$.

\section{The Roelcke compactification of groups of isometries}
\label{s:roelcomp}

Let $(M,d)$ be a complete 
$\o$-homogeneous Urysohn metric space, and let $G=\Is(M)$.
In the next section we shall
prove that $G$ is minimal and topologically simple.
The idea of the proof is to explicitly describe  the Roelcke compactification of
$G$. It turns out that the Roelcke completion of $G$ can be identified with
the compact space $\Th$ of all bi-Kat\v etov functions on $M^2$. 

In the preceding section we defined the embedding $i:G\to \Th$ by 
$i(\f)(x,y)=d(x,\f(y))$.
The space $\Th$, being compact,
has a unique compatible uniformity. Let $\sU$ be the coarsest uniformity on
$G$ which makes the map $i:G\to \Th$ uniformly continuous. We say that $\sU$
is the uniformity induced by $i$. The uniform space $(G,\sU)$ is isomorphic
to $i(G)$, considered as a uniform subspace of $\Th$. We are going to prove 
that $\sU$ is the Roelcke uniformity on $G$ (Theorem~\ref{th:Roel}).

Let us explain the idea of the proof. Let $\f, \f'\in G$. 
We want to prove that $\f$ and $\f'$ are
``sufficiently close" in $\Theta$ if and only if $\f'\in U\f U$, where
$U$ is a ``small" neighbourhood of the unity. Thus we are led to
the following question: under what conditions does the equation $\f'=\psi_1\f\psi_2$
have a solution with ``small" $\psi_1$ and $\psi_2$? Here ``small" means that points
of a given finite subset $A\sbs M$ are moved by less than $\e$. Observe that similar
questions for the equations $\f'=\f\psi$ or $\f'=\psi\f$ have an obvious answer:
$\f'\in\f U$ iff $\f$ and $\f'$ move points of $A$ ``almost in the same way", that
is, $d(\f(x),\f'(x))<\e$ for every $x\in A$; similarly, $\f'\in U\f$ iff the inverse maps
$\f\obr$ and ${\f'}\obr$ move points of $A$ ``almost in the same way".
The equation $\f'=\psi_1\f\psi_2$ with two unknowns $\psi_1$ and $\psi_2$
looks more complicated. However, the answer to the above question is easy also in this
case: the condition $\f'\in U\f U$ means that the finite metric spaces
$A\cup\f(A)$ and $A\cup\f'(A)$ are close to each other in the Gromov--Hausdorff metric.

We shall need the notion of the Gromov--Hausdorff metric only for finite metric
spaces with a given enumeration
(it differs from the usual notion dealing with non-enumerated spaces).
Let $X=\{x_1,\dots,x_n\}$ and $Y=\{y_1,\dots,y_n\}$
be two such spaces. The {\em Gromov--Hausdorff distance for enumerated spaces} 
between $X$ and $Y$, denoted
by $d_{GH}^{en}(X,Y)$, is the infimum of the numbers 
$\max\{D(x_i,y_i):i=1,\dots,n\}$, taken over all pseudometrics $D$ on $X\cup Y$ (we assume that 
$X$ and $Y$ are disjoint) such that $D$ induces the given metrics on $X$ and $Y$.
If $X$ and $Y$ have diameter $\le1$, we may assume that the same is true for $(X\cup Y, D)$,
otherwise replace $D$ by $D\wedge1$. Since the Urysohn space $(M,d)$ contains an isometric
copy of every finite metric space of diameter $\le1$ (Proposition~\ref{p:newequiv}), it follows 
that $d_{GH}^{en}(X,Y)$ is the infimum of the numbers $\max\{d(a_i,b_i):i=1,\dots,n\}$,
where $a_i, b_i\in M\ (1\le i\le n)$ are such that the correspondences $x_i\mapsto a_i$
and $y_i\mapsto b_i$ are isometric embeddings of $X$ and $Y$ into $M$, respectively.

\begin{prop}
\label{p:newGH}
Let $(X,d_X)$ and $(Y,d_Y)$ be two enumerated finite metric spaces, 
$X=\{x_1,\dots,x_n\}$, $Y=\{y_1,\dots,y_n\}$. Let 
$$
\e=\max\{|d_X(x_i, x_j)-d_Y(y_i,y_j)|: i,j=1,\dots,n\}.
$$
Then $d_{GH}^{en}(X,Y)=\e/2$
\end{prop}

\begin{proof}
The inequality $\ge$ is obvious: if $D$ is a pseudometric on $X\cup Y$ extending $d_X$ and $d_Y$
and $\e=|d_X(x_i, x_j)-d_Y(y_i,y_j)|$, 
then at least one of the numbers $D(x_i,y_i)$ and $D(x_j,y_j)$ must
be $\ge\e/2$. To prove the reverse inequality, we construct a pseudometric $D$ on $Z=X\cup Y$ extending $d_X$
and $d_Y$ such that
$$
D(x_i,y_i)=\e/2,\quad i=1,\dots,n.
$$
The function $D$ is defined by these requirements on $X^2$, $Y^2$, and the set $\{(x_i,y_i):i=1,\dots,n\}$.
To see that $D$ can be extended to a pseudometric on $Z$, it suffice to verify
that for any sequence $z_1,\dots,z_s$ of points of $Z$ such that
all the expressions $D(z_i,z_{i+1})$ ($1\le i <s$) and $D(z_1,z_s)$ are defined
the inequality
$$
D(z_1,z_s)\le \sum_{i=1}^{s-1} D(z_i,z_{i+1}) \leqno{\rm(A)}
$$ 
holds. Then the required extension is given by the formula
$$
D(z,z')=\inf\sum_{i=1}^{s-1}D(z_i,z_{i+1}),
$$
where the infimum is taken over
all chains $z_1=z,z_2,\dots,z_s=z'$ such that all the terms $D(z_i,z_{i+1})$
are defined.
An easy argument using induction shows that (A) follows from its special case:
for any ``quadrangle" in $Z$ of the form $x_i,\, y_i,\, y_j,\, x_j$ each of 
the four numbers
$d_X(x_i,x_j)$, $D(x_i, y_i)$, $d_Y(y_i,y_j)$, and $D(x_j,y_j)$ 
does not exceed
the sum of the three others. This case is obvious: for example,
since $d_X(x_i,x_j)-d_Y(y_i,y_j)\le\e$, we have
$$
d_X(x_i,x_j)\le d_Y(y_i,y_j)+\e=D(x_i, y_i)+d_Y(y_i,y_j)+D(x_j,y_j).
$$
%
%
\end{proof}

\begin{corol}
\label{c:GH}
Let $(X,d)$ be an Urysohn metric space. Let $a_1,\dots,a_n,b_1,\dots,b_n\in X$, and suppose that
$$
|d(a_i,a_j)-d(b_i,b_j)|\le 2\e
$$
for all $i,j=1,\dots,n$. Then there exist points $c_1,\dots,c_n\in X$ such that $d(c_i,c_j)=d(b_i,b_j)$
and $d(a_i,c_i)\le \e$ for all $i,j=1,\dots,n$.
\qed
\end{corol}

We now are in a position to prove the main result of this section.
Recall that $(M,d)$ is a complete 
$\o$-homogeneous Urysohn metric space, $G=\Is(M)$, and $\Theta$ is the space of 
bi-Kat\v etov functions on $M^2$ considered in the previous section.

\begin{thm}
\label{th:Roel}
The range of the embedding $i:G\to \Th$ is dense in $\Th$. The uniformity
$\sU$ on $G$ induced by the embedding $i$ coincides with the Roelcke uniformity
$\sL\wedge\sR$. Therefore,
$G$ is Roelcke-precompact, and the Roelcke compactification of $G$
can be identified with $\Th$.
\end{thm}

\begin{proof} 
If $A$ is a finite subset of $M$ and $\e>0$, 
let $U_{A,\e}=
\{\psi\in G: d(\psi(x),x)<\e\text{ for every } x\in A\}\in \sN(G)$.
Let $W_{A,\e}$ be the set of
all pairs $(f,g)\in \Th^2$ such that $|f(x,y)-g(x,y)|<\e$ for all $x,y\in A$.
The sets of the form $W_{A,\e}$ constitute a base of entourages of
the uniformity on $\Th$. If $(f,g)\in W=W_{A,\e}$, we say that $f$ and $g$ are
{\it $W$-close}.
Our proof proceeds in three parts.

(a) We prove that $i(G)$ is dense in $\Th$. Let $f\in \Th$, and let $Of$ be a
neighbourhood of $f$ in $\Th$. We must prove that $i(\f)\in Of$ for some
$\f\in G$. 

We may assume that $Of$ is the set of all $g\in \Theta$  such that
$g$ is $W_{A,\e}$-close to $f$:
$$
Of=\{g\in \Th: |g(x,y)-f(x,y)|<\e\text{ for all } x, y\in A\},
$$
where $A$ is a finite subset of $M$ and $\e>0$. Let $A=\{a_1,\dots, a_n\}$.
We claim that there exist points $b_1,\dots, b_n\in M$ such that
$d(b_i, b_j)=d(a_i, a_j)$ and $d(a_i, b_j)=f(a_i, a_j)$, $1\le i,j \le n$.
Indeed, since $f$ is bi-Kat\v etov, 
the formulas
above define a pseudometric on the set $F=\{a_1,\dots,a_n,b_1,\dots,b_n\}$,
where $b_1,\dots,b_n$ are new points. Since $M$ is Urysohn, the embedding of 
$A$ into $M$ extends to a distance-preserving map from $F$ to $M$.

Since $M$ is $\o$-homogeneous, there exists an isometry $\f$ of $M$ such
that $\f(a_i)=b_i$, $1\le i \le n$. Let $g=i(\f)$. For every $i,j\in [1,n]$
we have $g(a_i, a_j)=d(a_i, \f(a_j))=d(a_i, b_j)=f(a_i, a_j)$. Thus $g\in Of$.
This proves that $i(G)$ is dense in $\Th$.

(b)
We prove that the uniformity $\sU$ is coarser than $\sL\wedge\sR$. 

Whenever a topological group $H$ acts continuously 
on a compact space $X$
(on the left), for every $x\in X$ the orbit map $h\mapsto hx$ 
from $H$ to $X$ is 
right-uniformly continuous. We saw that $G$ acts continuously on $\Theta$
(Proposition~\ref{p:contshift}). The embedding $i:G\to \Theta$ can be viewed
as the orbit map corresponding to $d$, the neutral element of $\Theta$. 
It follows that $i$ is $\sR$-uniformly continuous. Similarly,
$i$ is $\sL$-uniformly continuous (use the right action of $G$ on $\Theta$,
or, alternatively, use the involution on $\Theta$ to deduce $\sL$-uniform
continuity of $i$ from its $\sR$-uniform continuity). 
Therefore, the uniformity $\sU$ is coarser than both $\sL$ and $\sR$
and hence coarser than $\sL\wedge\sR$. 

%

(c)
We prove that $\sU$ is finer than $\sL\wedge\sR$. It suffices to show that
for every $U\in \sN(G)$ 
there exists an entourage
$W$ of the uniformity on $\Th$ (in other words, a neighbourhood of the diagonal
of $\Th^2$) with the following property: if $\f,\f'\in G$ are such that 
$i(\f)$ and $i(\f')$ are $W$-close, then 
$\f'\in U\f U$. Assume that $U=U_{A,\e}$.
We claim that $W=W_{A,2\e}$ has the required property.

Let $\f,\f'\in G$ be such that $i(\f)$ and $i(\f')$ are $W_{A,2\e}$-close.
This means that 
$$
\d=\max\{|d(x,\f(y))-d(x,\f'(y))|: x,y\in A\}<2\e.
$$
Let $A=\{a_1,\dots, a_n\}$, $b_i=\f(a_i)$ and $c_i=\f'(a_i)$, $i=1,\dots, n$.
We have 
$d(b_i,b_j)=d(a_i, a_j)=d(c_i,c_j)$ and
$|d(a_i,b_j)-d(a_i,c_j|\le\d$ for all $i$ and $j$. 
In virtue of Corollary~\ref{c:GH}, there exist points $a_1', \dots, a_n',
b_1',\dots, b_n'\in M$ such that the correspondence $a_i\mapsto a_i'$,
$b_i\mapsto b_i'$ is distance-preserving and $d(a_i', a_i)\le \d/2<\e$, 
$d(b_i',c_i)\le \d/2<\e$. 


Since $M$ is
$\o$-homogeneous, there exists an isometry $\psi_1$ of $M$ such that
$\psi_1(a_i)=a_i'$ and $\psi_1(b_i)=b_i'$, $i=1,\dots,n$.
We have $\psi_1\in U$, since each $a_i$ is moved by less than $\e$.
Put $\psi_2=\f\obr\psi_1\obr\f'$. For every $i=1,\dots,n$ we have 
$d(\psi_2(a_i), a_i)=d(\f'(a_i),\psi_1\f(a_i))=d(c_i,b_i')<\e$, 
hence $\psi_2\in U=U_{A,\e}$. Thus $\f'=\psi_1\f\psi_2\in U\f U$, 
as required.
\end{proof}


Recall that a non-empty collection $\sF$ of non-empty subsets of a set $X$
is a {\it filter base\/} on $X$ if for every $A,B\in \sF$ there is $C\in \sF$
such that $C\sbs A\cap B$. If $X$ is a topological space, $\sF$ is a filter
base on $X$ and $x\in X$, then $x$ is a {\it cluster point\/} of $\sF$ if
every neighbourhood of $x$ meets every member of $\sF$, and $\sF$ 
{\it converges\/} to $x$ if every neighbourhood of $x$ contains a member of 
$\sF$. If $\sF$ and $\sG$ are two filter bases on $G$, let
$\sF\sG=\{AB:A\in \sF,\, B\in \sG\}$. 

For every $p\in \Th$ let $\sF_p=\{G\cap V: V\text{ is a neighbourhood of $p$ in
$\Th$}\}$. In other words, $\sF_p$ is the trace on $G$ of the filter of
neighbourhoods of $p$ in $\Th$.
If $p,q\in \Th$, it is not true in general that $\sF_p\sF_q$ converges to 
$p\bull q$. However, we have the following result, which will be used in the
proof of Theorem~\ref{th:1.5}:

\begin{prop}
\label{p:cluster}
If $p,q\in \Th$, then $p\bull q$ is a cluster point of the filter base
$\sF_p\sF_q$.
\end{prop}

\begin{proof}
Let $U_1$, $U_2$, $U_3$ be neighbourhoods of $p$, $q$ and $p\bull q$,
respectively. We must show that $U_3$ meets the set 
$(U_1\cap i(G))(U_2\cap i(G))$.

We may assume that for some finite set $A=\{a_1,\dots,a_n\}\sbs M$ 
and $\e>0$ we have
\begin{gather*}
U_1=\{f\in \Th: |f(x,y)-p(x,y)|<\e \text{ for all } x,y\in A\};\\
U_2=\{f\in \Th: |f(x,y)-q(x,y)|<\e \text{ for all } x,y\in A\};\\
U_3=\{f\in \Th: |f(x,y)-p\bull q(x,y)|<\e \text{ for all } x,y\in A\}.
\end{gather*}
We saw in the last paragraph of the preceding section
that there exist a metric space $(L,\rho)$ and isometric embeddings
$h_i:M\to L$ ($i=1,2,3$) such that $p(x,y)=\rho(h_1(x),h_2(y))$,
$q(x,y)=\rho(h_2(x),h_3(y))$ and $p\bull q(x,y)=\rho(h_1(x),h_3(y))$ for all
$x,y\in M$. Let $X=h_1(A)\cup h_2(A)\cup h_3(A)$. Since $M$ is Urysohn,
there exists an isometric embedding of $X$ into $M$ which extends the isometry
$h_1\obr:h_1(A)\to A$. It follows that there exist points 
$b_1,\dots,b_n,c_1\dots, c_n\in M$ such that 
$d(b_i,b_j)=d(c_i,c_j)=d(a_i,a_j)$, $d(a_i,b_j)=p(a_i,a_j)$, 
$d(b_i,c_j)=q(a_i,a_j)$ and $d(a_i,c_j)=p\bull q(a_i,a_j)$ for all $i,j$.
Since $M$ is $\o$-homogeneous, there exists an isometry $\f\in G$ such that
$\f(a_i)=b_i$, $1\le i \le n$. Let $x_i=\f\obr(c_i)$. Using again the
$\o$-homogeneity of $M$, we find an isometry $\psi\in G$ such that
$\psi(a_i)=x_i$, $1\le i\le n$. Note that $\f\psi(a_i)=c_i$ and $d(a_i,x_j)=
d(\f(a_i),\f(x_j))=d(b_i,c_j)=q(a_i,a_j)$ for all $i,j$. We claim that 
$i(\f)\in U_1$, $i(\psi)\in U_2$ and $i(\f\psi)\in U_3$. Indeed, we have
$i(\f)(x,y)=d(x,\f(y))=p(x,y)$ for all $x,y\in A$ and hence $i(\f)\in U_1$.
The other two cases are considered similarly. Thus
$i(\f\psi)\in \left((U_1\cap i(G))(U_2\cap i(G))\right)\cap U_3\ne\emptyset$.
\end{proof}

If $H$ is a group and $g\in H$, we denote by $l_g$ (respectively, $r_g$)
the left shift of $H$ defined by $l_g(h)=gh$ (respectively,
the right shift defined by $r_g(h)=hg$). 

\begin{prop}
\label{p:shiftroel}
Let $H$ be a topological group, and let $K$ be the Roelcke completion of $H$.
Let $g\in H$. 
Each of the following self-maps of $H$ extends to a self-homeomorphism of $K$:
(1)~the left shift~$l_g$; (2)~the right shift~$r_g$; (3)~the inversion 
$g\mapsto g\obr$; (4) the inner automorphism $h\mapsto ghg\obr$.
\end{prop}

\begin{proof}
Let $\sL$ and $\sR$ be the left and the right uniformity on $H$, respectively.
In each of the cases (1)--(4) 
the map $f:H\to H$ under consideration is an automorphism of
the uniform space $(H,\sL\wedge\sR)$. This is obvious for the cases~(3)
and~(4). For the cases~(1) and~(2), observe that the uniformities $\sL$ and
$\sR$ are invariant under left and right shifts, hence the same is true for
their greatest lower bound $\sL\wedge\sR$. It follows that in all cases
$f$ extends to an automorphism of the completion $K$ of the uniform space
$(H,\sL\wedge\sR)$.
\end{proof}

For the group $G$ and its Roelcke completion $\Th$ the validity
of Proposition~\ref{p:shiftroel} can be seen directly.
Recall that the embedding $i:G\to \Th$ is a morphism of monoids with 
an involution (see the two paragraphs before Proposition~\ref{p:contshift}).
The involution $f\mapsto f^*$ on $\Th$ is continuous and
hence coincides with the extension of the inversion on $G$ given by 
Proposition~\ref{p:shiftroel}. For every $g\in G$ let $L_g$, $R_g$ and $\Inn_g$
be the self-maps of $\Th$ defined by $L_g(p)=g\bull p$, $R_g(p)=p\bull g$ and
$\Inn_g(p)=g\bull p\bull g\obr$. 
These maps are extensions 
over $\Th$ of the left shift $l_g$ of $G$, the right shift $r_g$, and the inner
automorphism $l_g\circ r_{g\obr}$, respectively.
In virtue of Proposition~\ref{p:contshift}, the maps $L_g$ and $R_g$
are continuous, and the same is true for $\Inn_g=L_g\circ R_{g\obr}$.
%
%

An {\it inner automorphism\/} of $\Th$ is a map of the form $\Inn_g$, $g\in G$.
Proposition~\ref{p:contshift} shows that 
$\Inn_g(p)(x,y)=p(g\obr(x),g\obr(y))$ for all $p\in \Th$ and $x,y\in M$.
It follows that for every closed $F\sbs M$ we have $\Inn_g(b_F)=b_{g(F)}$,
where $b_F$ is the idempotent corresponding to $F$
(see Proposition~\ref{p:idemp}).

\begin{prop}
\label{p:innidemp}
There are precisely two idempotents in $\Th$ which are $\ge d$ and are
invariant under all inner automorphisms: the unity $d$ and the constant 1.
\end{prop}

\begin{proof}
According to Proposition~\ref{p:idemp}, every idempotent $\ge d$ is of the form
$b_F$ for some closed $F\sbs M$. If $b_F$ is invariant under 
inner automophisms, then $b_{g(F)}=\Inn_g(b_F)=b_F$ and hence $g(F)=F$ for
every $g\in G$. Since the action of $G$ on $M$ is transitive, no proper
non-empty subset of $M$ is $G$-invariant. Thus either $F=M$ or $F=\emptyset$.
Accordingly, either $b_F=d$ or $b_F=1$.   
\end{proof}

\section
{Proof of Theorem~\ref{th:1.5}}
\label{s:proof1.5}

We preserve the notation of the preceding section: $M$ is a complete
$\o$-homogeneous Urysohn metric space, $G=\Is(M)$, $\Th$ is the set of all
bi-Kat\v etov functions on $M^2$. We saw that $G$ is Roelcke-precompact and
that $\Th$ can be identified with the Roelcke compactification of $G$
(Theorem~\ref{th:Roel}). 
In this section we prove that $G$ is minimal and topologically simple.

\begin{prop}
\label{p:simplemin}
For every topological group $H$ the following
conditions are equivalent:
\begin{enumerate}
\item $H$ is minimal and topologically simple;
\item if $f:H\to H'$ is a continuous onto homomorphism of topological groups,
then either $f$ is a homeomorphism or $|H'|=1$.\qed
\end{enumerate}
\end{prop}

\begin{prop}
\label{p:compnorm}
The group $G$ has no compact normal subgroups other than $\{e\}$.
\end{prop}

We shall prove later that actually $G$ has no non-trivial closed normal 
subgroups.

\begin{proof}
Let $H\ne\{e\}$ be a normal subgroup of $G$. We show that $H$ is not compact.

Fix $a\in M$ and $f\in H$ such that $f(a)\ne a$. Let $r=d(f(a),a)$, and let
$S=\{x\in M:d(x,a)=r\}$ be the sphere of radius $r$ centered at $a$. We claim
that the orbit $Ha$ contains $S$. Fix $x\in S$. Since $M$ is $\o$-homogeneous,
there exists an isometry $g\in G$ which leaves the point $a$ fixed and maps
$f(a)$ to $x$. Let $h=gfg\obr$. Since $H$ is normal, we have $h\in H$ and hence
$x=h(a)\in Ha$. Thus $S\sbs Ha$, as claimed.

Since $M$ is Urysohn, we can construct by induction an infinite sequence
$x_1,x_2,\dots$ of points in $S$ such that all the pairwise distances between
distinct members of this sequence are equal to $r$. Since $S\sbs Ha$,
it follows that $Ha$ is not compact. Hence $H$ is not compact.
\end{proof}

Let $(L,\rho)$ be a metric space. A self-map $f:L\to L$ is 
{\it non-expanding\/} if \break
$\rho(f(x),f(y))\le \rho(x,y)$ for all $x,y\in L$.

\begin{lemma}
\label{l:nonexp}
Let $(L,\rho)$ be a metric space, and let
$F$ be the semigroup of all non-expanding self-maps of $L$, equipped 
with the topology of pointwise convergence. Then 
the map $(f,g)\mapsto f\circ g$ from $F^2$ to $F$ is continuous.
Thus $F$ is a topological semigroup.
\end{lemma}

This lemma and Proposition~\ref{p:Rajk} below are well known. We include
a proof for the reader's convenience.

\begin{proof}
It suffices to show that for every
$x\in L$ the map $(f,g)\mapsto f(g(x))$ from $F^2$ to $L$ is continuous. Fix
$f_0,g_0\in F$, $x\in L$ and $\e>0$. Let $y=g_0(x)$, 
$Of_0=\{f\in F:\rho(f(y),f_0(y))<\e\}$ and 
$Og_0=\{g\in F:\rho(g(x),y)<\e\}$. If $f\in Of_0$ and $g\in Og_0$, then
$\rho(f(g(x)),f_0(g_0(x)))\le \rho(f(g(x)),f(y))+\rho(f(y),f_0(y))
<\rho(g(x),y)+\e<2\e$. 
\end{proof}

\begin{prop}
\label{p:Rajk}
If $L$ is a complete metric space, then the group $\Is(L)$ is complete.
\end{prop}

Recall that we call a topological group complete if it is complete with respect
to the upper uniformity.

\begin{proof}
Let $X=L^L$ be the set of all self-maps of $L$, equipped with the product
uniformity. The group $H=\Is(L)$ can be considered as a subset of $X$.
The uniformity $\sU$ on $H$ induced by the product uniformity on $X$
coincides with the
left uniformity $\sL$. Indeed, a basic entourage for $\sU$ has the form
$W_{A,\e}=\{(f,g)\in H^2: \rho(f(x),g(x))<\e\text{ for all }x\in A\}$, 
where $\rho$ is the metric on $L$, $A$ is a finite subset of $L$ and $\e>0$. 
Let $U_{A,\e}=\{f\in H:\rho(f(x),x)<\e\text{ for all }x\in A\}$. Then
$U_{A,\e}$ is a basic neighbourhood of unity in $H$, 
and $W_{A,\e}=\{(f,g)\in H^2:g\obr f\in U_{A,\e}\}$ 
is a basic entourage for $\sL$. Thus $\sU=\sL$.
It follows that the map $g\to g\obr$ from $H$ to $X$ induces the 
right uniformity on $H$, and the map $j:H\to X^2$ defined by $j(g)=(g,g\obr)$
induces the upper uniformity $\sL\vee\sR$. Since $X^2$ is complete, to prove
that $H$ is complete it suffices to show that $j(H)$ is closed in $X^2$. 
Let $F$ be the set of all non-expanding self-maps of $L$. Then $F$ is
closed in $X$. The map $(f,g)\mapsto f\circ g$ from
$F^2$ to $F$ is continuous (Lemma~\ref{l:nonexp}). 
Since $j(G)=\{(f,g)\in F^2:fg=gf=\text{id}_L\}$,
it follows that $j(G)$ is closed in $F^2$ and hence in $X^2$.
\end{proof}

We say that a metric space $L$ is {\it homogeneous\/} if every point of $L$
can be mapped to every other point by an isometry of $L$ onto itself.

\begin{lemma}
\label{l:grweight}
If $L$ is a homogeneous metric space, then $w(\Is(L))=w(L)$.
\end{lemma}

\begin{proof}
For every metric space $X$ we have $w(\Is(X))\le w(X)$. If $X$ is homogeneous,
then for every $a\in X$ the map $f\to f(a)$ from $\Is(X)$ to $X$ is onto,
whence $w(X)\le w(\Is(X))$.
\end{proof}

We are now ready to prove Theorem~\ref{th:1.5}:

\smallskip
{\it
If $M$ is a complete $\o$-homogeneous
Urysohn metric space, then the group $G=\Is(M)$ is complete, 
Roelcke-precompact,
minimal and topologically simple. The weight of $G$ is equal to the
weight of $M$.
}
\smallskip

\begin{proof}
We saw that $G$ is Roelcke-precompact
(Theorem~\ref{th:Roel}). Proposition~\ref{p:Rajk} shows that
$G$ is complete, and Lemma~\ref{l:grweight} shows that $w(G)=w(M)$.
Let $f:G\to G'$ be a continuous onto homomorphism.
According to Proposition~\ref{p:simplemin}, to prove that $G$ is minimal and
topologically simple, it suffices to prove that either $f$ is a homeomorphism
or $|G'|=1$. 

Since $G$ is Roelcke-precompact, so is $G'$. Let $\Th'$ be the Roelcke
compactification of $G'$. The homomorphism $f$ extends to a continuous map
$F:\Th\to\Th'$. Let $e'$ be the unity of $G'$, and let $S=F\obr(e')\sbs \Th$.

\smallskip 

{\em Claim 1.} $S$ is a subsemigroup of $\Th$.

\smallskip 

Let $p,q\in S$. In virtue of Proposition~\ref{p:cluster}, there exist filter
bases $\sF_p$ and $\sF_q$ on $G$ such that $\sF_p$ converges to $p$ (in $\Th$),
$\sF_q$ converges to $q$ and $p\bull q$ is a cluster point of the filter base
$\sF_p\sF_q$. The filter bases $\sF_p'=F(\sF_p)$ and $\sF_q'=F(\sF_q)$ on $G'$
converge to $F(p)=F(q)=e'$, hence the same is true for the filter base
$\sF_p'\sF_q'=F(\sF_p\sF_q)$. Since $p\bull q$ is a cluster point of
$\sF_p\sF_q$, $F(p\bull q)$ is a cluster point of the convergent filter base
$F(\sF_p\sF_q)$. A convergent filter on a Hausdorff space has only one cluster
point, namely the limit. Thus $F(p\bull q)=e'$ and hence $p\bull q\in S$.

\smallskip 

{\em Claim 2.} The semigroup $S$ is closed under involution.

\smallskip 

In virtue of Proposition~\ref{p:shiftroel},
the inversion on $G'$ extends to an involution $x\mapsto x^*$ of $\Th'$. Since
$F(p^*)=F(p)^*$ for every $p\in G$, the same holds for every $p\in \Th$. Let
$p\in S$. Then $F(p^*)=F(p)^*=e'$ and hence $p^*\in S$.

\smallskip 

{\em Claim 3.} If $g\in G$ and $g'=f(g)$, then $F\obr(g')=g\bull S=S\bull g$.

\smallskip 

We saw that the left shift $h\mapsto gh$ of $G$ extends to a continuous
self-map $L=L_g$ of $\Th$ defined by $l(p)=g\bull p$
(Proposition~\ref{p:contshift}). 
According to Proposition~\ref{p:shiftroel}, the self-map $x\mapsto g'x$ of $G'$
extends to a self-homeomorphism $L'$ of $\Th'$. The maps 
$F\circ L$ and $L'\circ F$ from
$\Th$ to $\Th'$ coincide on $G$ and hence everywhere. Replacing
$g$ by $g\obr$, we see that $F\circ L\obr=(L')\obr\circ F$.
Thus
$F\obr(g')=F\obr L'(e')=LF\obr(e')=g\bull S$. Using right shifts,
we similarly conclude that $F\obr(g')=S\bull g$.

\smallskip 

{\em Claim 4.} $S$ is invariant under inner automorphisms of $\Th$.

\smallskip 

We have just seen that $g\bull S=S\bull g$ for every $g\in G$, hence
$g\bull S\bull g\obr=S$.

\smallskip 

Let $T=\{f\in S:f\ge d\}$. 
Note that $i(e_G)=d\in T\ne\emptyset$.
According to Proposition~\ref{p:maximidem},
there is a greatest element $p$ in $T$, and $p$ is idempotent. 
Since inner automorphisms of $\Th$ preserve the order on $\Th$ 
and the unity $d$, Claim 4
implies that $p$ is invariant under inner automorphisms.
In virtue of Proposition~\ref{p:innidemp}, either $p=d$ or $p=1$.
We shall show that either $f$ is a homeomorphism or $|G'|=1$,
according to which of the cases $p=d$ or $p=1$ holds.

Consider first the case $p=d$.

\smallskip 

{\em Claim 5.} If $p=d$, then all elements of $S$ are invertible in $\Th$.

\smallskip 

Let $f\in S$. Then $f^*\bull f\in S$
and $f\bull f^*\in S$, since $S$ is a symmetrical semigroup.
According to Proposition~\ref{p:charbikat}, we have 
$f^*\bull f\ge d$ and $f\bull f^*\ge d$.
Since $p=d$, there are no elements $>d$ in $S$. Thus the inequalities
$f^*\bull f\ge d$ and $f\bull f^*\ge d$ are actually equalities. It follows
that $f^*$ is the inverse of $f$.

\smallskip 

{\em Claim 6.} If $p=d$, then $S=\{e\}$.

\smallskip 

Claim 5 and Proposition~\ref{p:invert} imply that $S$ is a subgroup of $G$.
This subgroup is normal (Claim~4) and compact, since $S$ is closed in $\Th$.
Proposition~\ref{p:compnorm} implies that $S=\{e\}$.

\smallskip 

{\em Claim 7.} If $p=d$, then $f:G\to G'$ is a homeomorphism.

\smallskip 

Claims~6 and~3 imply that $G= F\obr(G')$ and that the map $f:G\to G'$ is
bijective. Since $F$ is a map between compact spaces, it is perfect, and hence
so is the map $f:G=F\obr(G')\to G'$. Thus $f$, being a perfect bijection,
is a homeomorphism.

\smallskip 

Now consider the case $p=1$.

\smallskip 

{\em Claim 8.} If $1\in S$, then $G'=\{e'\}$.

\smallskip 

Let $g\in G$ and $g'=f(g)$. We have $g\bull 1=1\in S$. On the other hand,
Claim~3 implies that $g\bull 1\in g\bull S=F\obr(g')$. Thus $g'=F(g\bull 1)=
F(1)=e'$.
\end{proof}

\section
{Remarks}
\label{s:rem}

1.
Let $M$ be a complete $\o$-homogeneous Urysohn metric space, and let
$G=\Is(M)$. In Section~\ref{s:roelcomp}
we identified the Roelcke completion of $G$
with the set $\Th$ of all bi-Kat\v etov functions on $M^2$. The set $\Th$
was equipped with structures of three kinds: topology, order, semigroup
structure. The proof of Theorem~\ref{th:1.5} was based on the interplay between
these three structures. We now establish a natural one-to-one correspondence
between $\Th$ and a set of closed relations on a compact space. This
correspondence will be an isomorphism for all three structures on $\Th$.

Let $K$ be a compact space. A {\it closed relation\/} on $K$
is a closed subset of $K^2$.
Let $E(K)$ be the compact space  of all closed relations on $K$,
equipped with the Vietoris topology. The set $E(K)$ has a natural partial 
order. If $R,S\in E(K)$, then the composition
$R\circ S$ is a closed relation, since $R\circ S$ is the image of the closed
subset $\{(x,z,y): (x,z)\in S,\ (z,y)\in R\}$ of $K^3$ under the projection
$K^3\to K^2$ which is a closed map. Thus $E(K)$ is a semigroup
with involution. In general the map $(R,S)\mapsto R\circ S$ from
$E(K)^2$ to $E(K)$ is not separately continuous (neither left nor right
continuous).

We denote by $\Homeo(K)$
the group of all
self-homeomorphisms of $K$, equipped with the compact-open topology.
For  every  $h\in  \Homeo(K)$ let $\G(h)=\{(x, h(x)): x\in K\}$ be the graph of
$h$.
The map $h\mapsto \G(h)$ from $\Homeo(K)$ to $E(K)$ is a homeomorphic
embedding and a morphism of monoids with an involution. The uniformity
induced on $\Homeo(K)$ by this embbedding is coarser than the Roelcke uniformity.

Now let $K$ be the compact space of all non-expanding
functions  $f:M\to I=[0,1]$, considered as a subspace of the product
$I^M$. There is a natural left action of $G$ on $K$, defined by  $gf(x)=f(g\obr
(x))\  (g\in  G,\  f\in  K,  \  x\in  M)$.  This action gives rise to
a morphism $G\to \Homeo(K)$ of topological groups which is easily seen to be
a homeomorphic embedding. Let $j:G\to E(K)$ be the composition of
this embedding with the map $h\mapsto \G(h)$ from $\Homeo(K)$ to $E(K)$. If
$g\in G$, then $j(g)$ is the relation $\{(f,gf):f\in K\}$. Let $\Phi$ be the
closure of $j(G)$ in $E(K)$. Let $\Th$ and $i:G\to \Th$ be the same as in
Sections~6 and~7.

\begin{thm}
\label{th:Rel}
The uniformity on $G$ induced by the embedding $j:G\to E(K)$ coincides with
the Roelcke uniformity, hence $\Phi$ can be identified with the Roelcke
compactification of $G$. The set $\Phi$ is a subsemigroup of $E(K)$. There
exists a unique homeomorphism $H:\Phi\to \Th$ such that $i=Hj$. The map
$H$ is an isomorphism of ordered semigroups.
\end{thm}

We omit the detailed proof and confine ourselves by a description of
the isomorphism $H$. If $R\in \Phi$, let $H(R)$ be the bi-Kat\v etov function
on $M^2$ defined by $H(R)(x,y)=\sup\{|q(x)-p(y)|:(p,q)\in R\}$, $x,y\in M$.
If $f\in \Th$, the relation $H\obr(f)$ is defined by
$H\obr(f)=\{(p,q)\in K^2: |q(x)-p(y)|\le f(x,y)\text{ for all } x,y\in M^2\}$.

Let us see what some of the results about $\Th$ obtained in Section~6
mean in terms of relations on $K$. Functions $p\in \Th$ which are $\ge d$
correspond (via the isomorphism $H$) to relations $R\in \Phi$ which contain
the diagonal of $K^2$ or, in other words, are reflexive. Thus 
Proposition~\ref{p:charbikat} 
implies that for every $R\in \Phi$ the relations $RR\obr$ and
$R\obr R$ are reflexive. This is equivalent to the fact that for every 
$R\in \Phi$ the domain and the range of $R$ is equal to $K$.

According to Proposition~\ref{p:idemp}, each idempotent $\ge d$ in $\Th$
has the form $b_F$ for some closed $F\sbs M$. Note that each  $b_F$ 
is symmetrical. Symmetrical idempotents $\ge d$ in $\Th$ correspond to 
relations in $\Phi$ which are reflexive, symmetrical and transitive or, 
in other words, are equivalence relations. For every closed $F\sbs M$ let
$R_F=H(b_F)$ be the equivalence relation corresponding to the idempotent
$b_F$. Two non-expanding functions $f,g\in K$ are $R_F$-equivalent if and
only if $f|F=g|F$. Proposition~\ref{p:idemp} implies that an equivalence
relation $R$ on $K$ belongs to $\Phi$ if and only if $R=R_F$ for some closed
$F\sbs M$.

2. Let $H$ be a Hilbert space, and let $G=U(H)$ be the group of all unitary 
operators on $H$, equipped with the pointwise convergence topology.
L.~Stoyanov proved that
$G$ is totally minimal \cite{Sto, DPS}. The methods of the present paper
yield an alternative proof of this theorem. Let $\sB(H)$ be the algebra
of all bounded linear operators on $H$. 
The {\it weak operator topology\/}
on $\sB(H)$ is the coarsest topology  such that for every $x,y\in H$ 
the function $A\mapsto (Ax,y)$ on $\sB(H)$ is continuous. 
Let $T=\{A:\Vert A\Vert\le 1\}$
be the unit ball in $\sB(H)$, equipped with the weak operator topology.
The Roelcke compactification of the group $G$ 
can be identified with $T$ \cite{Orsat}.
The set $T$ is a semigroup, and the idempotents in $T$ are
the orthogonal projectors. 
The proof of the
fact that $G$ is totally minimal proceeds similarly to the proof of
Theorem~\ref{th:1.5}.  Let us indicate the main steps. 
Let $f:G\to G'$ be a surjective morphism of topological groups. 
To prove that $G$ is totally minimal, it suffices to prove that $f$ is a 
quotient map. Extend $f$
to a map $F:T\to T'$, where $T'$ is the Roelcke compactification of $G'$. 
Let $e'$ be the unity of $G'$, and
let $S=F\obr(e')$ be the kernel of $F$. Then $S$ is a closed subsemigroup of
$T$. It turns out that
every closed subsemigroup of $T$ contains a least idempotent.
Let $p$ be the least idempotent in $S$. Since $S$ is invariant under
inner automophisms of $\Th$, so is $p$. It follows that either $p=1$ or $p=0$.
If $p=1$, then
$S\sbs G$, $G=F\obr (G')$ and the map $f$ is
perfect. If $p=0$, then $G'=\{e'\}$. See \cite{Orsat} for more details.

3. Our method of proving minimality, based on the consideration of the
Roelcke compactifications, can be applied to some groups of homeomorphisms.
A zero-dimensional compact space $X$ is {\it $h$-homogeneous\/} if all 
non-empty clopen subsets of $X$ are homeomorphic to each other. Let $K$
be a zero-dimensional $h$-homogeneous compact space, and let $G=\Homeo(K)$.
Then $G$ is minimal and topologically simple \cite{U4}.
Let us sketch a proof of this fact
which closely follows the proof of Theorem~\ref{th:1.5}. In the special case
when $K=2^\o$ is the Cantor set, the minimality of $\Homeo(K)$
was proved by Gamarnik \cite{Gamarn}.

Let $T$ be the compact space of all closed relations $R$ on $K$ such that 
the domain and the range of $R$ are equal to $K$. The map
$h\mapsto \G(h)$ from $G$ to $T$ 
induces the Roelcke uniformity on $G$, and the range $\G(G)$ of this map
is dense in $T$. Thus the Roelcke
compactification of $G$ can be identified with $T$. 
We noted that the set $E(K)$ of all closed relations on $K$
is an ordered semigroup with an involution.
The set $T$ is a closed symmetrical subsemigroup of $E(K)$.
Let $\D$ be the diagonal in $K^2$.
A relation $R\in E(K)$ is an equivalence relation if
and only if $R$ is a symmetrical idempotent and $R\supset\D$.
Let $S$ be a closed subsemigroup of $E(K)$, and let $S_1$ be the set of all
$R\in S$ such that $R\supset \D$. The proof of Proposition~\ref{p:maximidem}
shows that the set $S_1$, if it is non-empty,
has a largest element $P$, and $P$ is an idempotent.
If $S$ is symmetrical, then so is $P$, hence $P$ is an equivalence relation.

Now let $f:G\to G'$ be a surjective morphism of topological groups.
We show that either $f$ is a homeomorphism or $|G'|=1$. 
According to Proposition~\ref{p:simplemin},
this means that
$G$ is minimal and topologically simple. Extend $f$
to a map $F:T\to T'$, where $T'$ is the Roelcke compactification of $G'$. 
Let $e'$ be the unity of $G'$, and
let $S=F\obr(e')$. Then $S$ is a closed symmetrical subsemigroup of
$T$. Let $P$ be the largest element in the set $S_1=\{R\in S:\D\sbs R\}$.
Then $P$ is an equivalence relation on $K$. Since $S$ is $G$-invariant, so is
$P$. But there are only two $G$-invariant closed
equivalence relations on $K$, namely
$\D$ and $K^2$. If $P=\D$, then $S\sbs G$, $G=F\obr(G')$ and $f$ is perfect.
Since $G$ has no non-trivial compact normal subgroups, we conclude that $f$
is a homeomorphism. If $P=K^2$, then $S=T$ and $G'=\{e'\}$.

It is not clear if a similar argument can be used when $K$ is
a Hilbert cube and $G=\Homeo(K)$, see Problem~\ref{pro:minim}.

\section
{Acknowledgement}
\label{s:ack}

I am much obliged to the referee for careful reading of the paper and
for suggesting quite a few improvements.

\end{document}